\documentclass
{article}

\usepackage{graphicx}%
\usepackage{multirow}%
\usepackage{amsmath,amssymb,amsfonts}%
\usepackage{amsthm}%
\usepackage{mathrsfs}%
\usepackage[title]{appendix}%
\usepackage{xcolor}%
\usepackage{textcomp}%
\usepackage{manyfoot}%
\usepackage{booktabs}%
\usepackage{algorithm}%
\usepackage{algorithmicx}%
\usepackage{algpseudocode}%
\usepackage{listings}%
\usepackage{bm}
\usepackage{caption}

\begin{document}


\newtheorem{Th}{Theorem}[section]
\newtheorem{Cor}{Corollary}[section]
\newtheorem{Prop}{Proposition}[section]
\newtheorem{Lem}{Lemma}[section]
\newcounter{Remark}[section]
\newcounter{Example}[section]


\newcommand{\bet}{\begin{Th}}
\newcommand{\ent}{\stepcounter{Cor}
   \stepcounter{Prop}\stepcounter{Lem}
   \stepcounter{Remark}\stepcounter{Example}\end{Th}}


\newcommand{\bec}{\begin{Cor}}
\newcommand{\enc}{\stepcounter{Th}
   \stepcounter{Prop}\stepcounter{Lem}
   \stepcounter{Remark}\stepcounter{Example}\end{Cor}}
\newcommand{\bep}{\begin{Prop}}
\newcommand{\enp}{\stepcounter{Th}
   \stepcounter{Cor}\stepcounter{Lem}
   \stepcounter{Remark}\stepcounter{Example}\end{Prop}}
\newcommand{\bel}{\begin{Lem}}
\newcommand{\enl}{\stepcounter{Th}
   \stepcounter{Cor}\stepcounter{Prop}
   \stepcounter{Remark}\stepcounter{Example}\end{Lem}}

\newcommand{\Remark}{
   \stepcounter{Remark}
   \noindent{\bf Remark\,  \thesection.\theRemark \, }:
   \stepcounter{Th}\stepcounter{Cor}\stepcounter{Prop}
   \stepcounter{Lem}\stepcounter{Example}}
\newcommand{\Example}{
   \stepcounter{Example}
   \noindent{\bf Example\,  \thesection.\theExample \, }:
   \stepcounter{Th}\stepcounter{Cor}\stepcounter{Prop}
   \stepcounter{Lem}\stepcounter{Remark}}
\newcommand{\Proof}{\noindent{\it Proof\,}:\ }

\newcommand{\arc}[1]{%
 \settowidth{\dimen0}{\ensuremath{#1}}%
 \divide\dimen0 by 2%
 \overset{\rotatebox{-90}{\ensuremath{\left(\rule{0pt}{\dimen0}\right.}}}{#1}%
}
\newcommand{\overarc}[1]{\stackrel{\Large\mbox{$\frown$}}{#1}}


\def\qed{\hfill $\Box$}
\def\R{\mathbb{R}}
\def\Z{\mathbb{Z}}
\def\C{\mathbb{C}}
\def\F{\mathbb{F}}
\def\H{\mathbb{H}}
\def\Q{\mathbb{Q}}
\def\eps{\epsilon}
\def\La{\Lambda}
\def\al{\alpha}
\def\be{\beta}
\def\ga{\gamma}
\def\de{\delta}
\def\si{\sigma}
\def\om{\omega}
\def\x{\bm{x}}
\def\n{\bm{n}}
\def\la{\lambda}
\def\bary{\begin{array}{ll}}
\def\eary{\end{array}}
\def\beq{\begin{equation}}
\def\enq{\end{equation}}
\def\bary{\begin{array}{ll}}
\def\eary{\end{array}}
\def\qed{\hfill $\Box$}
\def\la{\lambda}
\def\ve{\varepsilon}
\def\bpm{\begin{pmatrix}}
\def\epm{\end{pmatrix}}
\def\vp{\varphi}

\date{}

\title
{Chern's Conjecture in the  Dupin case}

\author{Reiko Miyaoka} 

\maketitle


\abstract{Chern's  conjecture states that a closed minimal hypersurface in the euclidean sphere is isoparametric 
if it has constant scalar curvature. 
When the number $g$ of distinct principal curvatures is greater than three, 
few satisfactory results have been known. 
We attack the conjecture in the Dupin hypersurface case.  Our results are: 
A closed proper Dupin hypersurface with constant mean curvature is isoparametric (i) if $g=3$, 
 (ii) if $g=4$ and has constant scalar curvature, or (iii) if $g=4$ and has constant Lie curvature, and   
 (iv) if $g=6$ and has constant Lie curvatures. 
These cover all the non-trivial cases for a closed proper Dupin  
to be isoparametric since $g$ can take only values $1,2,3,4,6$. 
The originality of the proof is a use of topology and geometry, 
 which reduces assumptions needed in the algebraic argument.}\footnote
{{\it 2010 Mathematics Subject Classification:} Primary: 53C40. Secondary: 53C42, 53A07. 

{\it Keywords:} Chern's conjecture, Dupin hypersurface, isoparametric hypersurface, Lie sphere geometry, Lie curvature

Supported by Grants-in-Aid for Scientific Research C 21K03214 and 24K06701  of JSPS.}

\section{Introduction}
\label{sintro}

Wave fronts developing at a constant speed are called isoparametric hypersurfaces.  
Typical examples are light waves, and the subject's origin is geometric optics in the early 20th in Italy \cite{L}, \cite{So}. 
In the Euclidean and hyperbolic spaces, they are only concentric spheres and cylinders \cite{C1}. 
\'E. Cartan first treated the problem uniformly in the space forms  
and found more examples in the Euclidean sphere $S^{n}$ \cite{C2}. 

Such hypersurfaces exist in a family of parallel constant mean curvature (CMC) hypersurfaces. 
In fact, they have constant {\em principal} curvatures if the ambient space is a real space form. 
Each hypersurface is a piece of a complete hypersurface, namely, 
extendable to a global wave front.  

In the case $S^{n}$, the number $g$ of distinct principal curvatures is limited to 1,2,3,4,6.  
They are completely classified (Yau's 34th problem \cite{Y2}) and well investigated 
\cite{C1}, \cite{C2}, \cite{CCJ}, \cite{Chi}, \cite{DN}, \cite{FKM},  \cite{Mdn}, \cite{Mg6}, \cite{Mg6E}, 
\cite{Mu1}, \cite{Mu2}. 

When a hypersurface  $M$ in $S^{n}$ has  constant principal curvatures $\la_1$,  $\la_2$, $\dots$,  $\la_{n-1}$,  
$M$ has constant scalar curvature $R_M=(n-1)(n-2)+H^2-\|A\|^2=(n-1)(n-2)+(\sum\la_i)^2-\sum\la_i^2$ 
where $A$ is the shape operator of $M$ and $H=$Tr$A$ is the mean curvature.  

\smallskip\noindent
{\bf Chern's conjecture} (\cite{Ch}): 
{\em A closed minimal hypersurface in the euclidean sphere is isoparametric 
if it has constant scalar curvature (CSC). }

\smallskip
This is Yau's 105th problem \cite{Y1} and is still open in general. 
Almeida-Brito affirmatively solves the first non-trivial case $\dim M=3$ in 
\cite{dAB} together with  Peng-Terng \cite{PeT}, Chang \cite{Cha} and Cheng-Wan \cite{CW}. 
When $\dim M>3$, Chen-Li give an affirmative answer even locally for $g\le3$  \cite{CL}. 
However, for $g\ge4$, no definite conclusion has been given without further assumptions \cite{SWY}.   

A hypersurface is called {\em Dupin} if
each principal curvature has constant multiplicity and 
 is constant along its curvature direction. 
Pinkall calls it a proper Dupin in \cite{P}, but here we omit ``proper".  
Instead, if a hypersurface satisfies only the second condition, we call it a weak Dupin. 

It is known \cite{Th} that a closed embedded Dupin hypersurface has 
the same topology as isoparametric hypersurfaces 
where ``closed" means a compact without boundary.  
In particular, $g\in\{1,2,3,4,6\}$ holds with some restricted multiplicities. 
 
In general, a Dupin hypersurface is not necessarily a piece of complete hypersurface. 
Moreover, they belong to the category of Lie sphere geometry, 
wider than Riemannian and conformal geometry. 
For details, see \cite{Ce}, \cite{P}, and \S\ref{sLiec}. 

 When $g=1,2,3$, it is known that Dupin hypersurfaces are images of isoparametric hypersurfaces 
under Lie contact transformations \cite{CR1},  \cite{CJ}. 
The last one is the extension of the author's result in the closed case \cite {M3},   
which motivated Cecil and Ryan to conjecture that closed embedded Dupin 
hypersurfaces are Lie images of isoparametric hypersurfaces \cite{CR1}. 
However,  we construct {\em counterexamples}  
for $g=4,6$ in  \cite{MO} (see also \cite{PT} for $g=4$). 
Here the invariant ``Lie curvature" discovered by the author \cite{M4} plays an essential role. 
The Lie curvature is given as the cross ratio of four distinct principal curvatures (\S\ref{sLiec}).  

Therefore, it seems still hard to attack Chern's conjecture in the Dupin case with more than three principal curvatures. 

We state our main results:  
 \bet  
A closed CMC Dupin hypersurface  $M$ in $S^{n}$ is isoparametric  
  \begin{enumerate}
\item[\textup {(i)}]  if $g=3$.  

 \item[\textup {(ii)}]  if $g=4$ and $M$  has  constant scalar curvature.  
 
\item[\textup{(iii)}]  if $g=4$ and $M$  has  constant Lie curvature.

\item[\textup{(iv)}]  if $g=6$ and $M$  has constant Lie curvatures. 

\end{enumerate} 
  \label{tmain}
  \ent

Due to (ii), 
Chern's conjecture is weakened to: A closed CMC  hypersurface with $g=4$  
is Dupin if it has CSC. 

In (iii) and (iv), we use the Lie curvatures instead of the scalar curvature.   
When $ g=4 $, the Lie curvature is essentially unique. 
When $g=6$, three Lie curvatures are essential if they are independent, 
i.e., if any one of the three cannot be expressed by the other two (Lemma in the front page of \cite{M6}). 
 Thus the condition in (iv)  can be replaced by ``if $g=6$ 
 and $M$ has three independent constant Lie curvatures". 

When $g=4,6$, even if we assume all the Lie curvatures are constant, it does not imply 
the Lie equivalence between isoparametric and Dupin hypersurfaces \cite{M4}, \cite{M6}.  
Thus the Dupin is still far from the isoparametric.   

 In the proof, the topological global property ``tautness" of closed Dupin hypersurfaces plays an important role.   
Note that  taut implies weak Dupin \cite{Th}, \cite{M2}, \cite{P2}. 
    
 \smallskip
On the other hand, for a hypersurface in $S^n$, 
 a principal curvature with multiplicity greater than one 
 is constant along its curvature direction \cite{O}, \cite{R}.  
 Thus we obtain: 
  \bec 
 Let $M$ be a closed CMC  hypersurface in $S^{n}$ with distinct $g$ principal curvatures 
 of which multiplicities are constant and greater than one. Then $M$ is isoparametric, 
 \begin{enumerate}
\item[\textup{(i)}]   if $g=3$.  
\item[\textup{(ii)}]  if $g=4$ and  $M$  has constant scalar curvature. 
\item[\textup{(iii)}]   if  $g=4$ and has constant  Lie curvature.
\item[\textup{(iv)}]   if  $g=6$ and  has three independent constant Lie curvatures. 
\end{enumerate} 
 \label{cmain}
 \enc

 For further results on Chern's conjecture, see \cite{SWY},  \cite{TY} and the references therein. 
  
  The paper is organized as follows:
  In \S\ref{spre}, we give a basic definition and the known results. 
  In \S\ref{sclosedg3} and in  \S\ref{sg4}, 
  we give proofs of (i) and (ii) of Theorem \ref{tmain}, which are rather elementary. 
However, the proofs of (iii) and (iv) need the Lie sphere geometry as an essential tool, 
and an advanced argument is necessary. 
In \S\ref{sLiec}, we briefly introduce the Lie sphere theory and describe the Lie curvature.   
In \S\ref{sg4iii} and in  \S\ref{sg6}, we prove (iii) and (iv) of Theorem \ref{tmain}.

\section{Preliminaries}
\label{spre}

We mainly follow the notation in \cite{M4}. 
Let $S^{n}$ be the $n$-dimensional unit sphere in $\R^{n+1}$ with center at the origin. 
Consider an isometrically embedded orientable hypersurface $\bm{p}: M\to  S^{n}$ with 
a unit normal vector field $n$. 
Let $A$ be the second fundamental tensor of  $M$ 
and $\la_1\ge\dots\ge\la_{n-1}$ be the principal curvatures. 
For $\la\in\{\la_i\}$,  the curvature distribution $D(\la)$ is defined to be 
\beq
D_p(\la)=\{X\in T_pM\mid AX=\la X\}, \quad p\in M.
\label{Distribution}
\enq
 We abuse $p$ for $\bm{p}(p)$ if it is clear. 
 The following is well-known \cite{O}, \cite{R}.

\smallskip\noindent
{\bf Fact 1.}  When $\dim D(\la)$ is constant ($= m$) on $M$, $D(\la)$ is an involutive distribution. 
Moreover if $\la$ is constant along $D(\la)$, (this is the case when $m> 1$), 
 the leaf $L$ is a piece of an $m$-dimensional subsphere of the curvature sphere $C(\la)$ at $p$.  

Recall we call a hypersurface {\em Dupin} if each principal curvature has constant multiplicity 
and is constant along its curvature direction. 

\Example The following are Dupin hypersurfaces:

(1) Isoparametric hypersurfaces. 

(2) Conformal or Lie images (\S\ref{sLiec}) of isoparametric hypersurfaces.   

(3) When $M$ is a Dupin hypersurface in $\R^m$,  
 take a cylinder or  a tube over $M$ in $\R^m\oplus \R^k$.  
Then we have a new (weak) Dupin hypersurface in $\R^{m+k}$.  
By a stereographic projection, we have a weak Dupin in $S^{m+k}$. 

\Remark (3) implies that {\em locally}, there are Dupin hypersurfaces with any number of principal curvatures with 
any  multiplicities  \cite{P}.

\smallskip
The focal points of $p\in M$ are given by 
\beq
f_p^i=\cos\theta_i\,p+\sin\theta_i\,n_p,\quad 0<\theta_i=\cot^{-1}\la_i<\pi
\label{cocal}
\enq
and by their antipodals $\bar f_p^i$  in $S^{n}$, $1\le i \le n-1$. 
For $x\in S^{n}$, let $l_x: M\to \R$ be the  function on $M$ defined   
by 
\beq
 l_x(p)=d(x,p)^2
 \label{distance}
 \enq
 where $d(x,p)$ is the spherical distance from $x$ to $p$. 
Let $p\in M$ be a critical point of $l_x$. 
Then the index of $p$ is given by the sum of the multiplicities of the focal points of $p$ 
on $\stackrel{\frown}{xp}$, where $\stackrel{\frown}{xp}$ is  
the oriented open geodesic segment from $x$ to $p$. 
For details, see \cite{CR1}. 
Remarkably, a closed embedded Dupin hypersurface is {\it taut} \cite{M2}, \cite{P2}, by which we mean
$l_x$ satisfies the equality in the Morse inequality for generic $x\in S^n$, where we use the homology 
with $\Z_2$ coefficient. 
In such a case, $l_x$ is called a {\it perfect Morse function}. 

It is obvious  to see  (locally): 

\smallskip\noindent
(1) {\em When $g=1$, a Dupin hypersurface is a piece of hypersphere.} 

\noindent 
(2) {\em When $g=2$, a CMC Dupin is a piece of a Clifford hypersurface} (orthogonal product of two spheres), 
i.e., of an isoparametric hypersurface with $g=2$. 
\label{p12}
 
 \smallskip
In fact,  when $g=2$,   CMC implies that  the principal curvatures $\la,\mu$ satisfy 
$m_1\la+ m_2\mu=H$ for some constant $H$, where $m_1$ ($m_2$, resp.) is the multiplicity of $\la$ ($\mu$, resp.).  
Since $\la$ ($\mu$, resp.) is constant along  $D(\la)$ ($D(\mu)$, resp.),  
 so is $\mu$ ($\la$, resp.), and hence $\la,\mu$ are constant all over $M$. 
 
\smallskip
When $\dim M=3$, Chern's conjecture is solved affirmatively 
in  \cite {dAB}, \cite{PT},  \cite{Ch}, and \cite{CW}. 
A local version for $g\le3$, see \cite{CL}. 
In Theorem \ref{tmain} (i), we treat all closed CMC hypersurface $M$ with $g=3$,   
where instead of constant scalar curvature, we assume Dupin. 
Our argument will help understand the cases $g=4,6$. 

\noindent
{\bf Fact 2.}
Isoparametric hypersurfaces in $S^n$ consist of  parallel CMC hypersurfaces $M_\theta$ 
with constant principal curvatures $\la_1>\dots>\la_g$,  $g\in\{1,2,3,4,6\}$ given by \cite{Mu1}
\beq
\la_i(\theta)=\cot \theta_i=\cot\Big(\dfrac{\pi}{2g}+\theta+\dfrac{(i-1)\pi}{g}\Big),\,\, i\in\{1,\dots,g\}, \,
-\dfrac{\pi}{2g}<\theta<\dfrac{\pi}{2g}. 
\label{thetai}
\enq
In particular, we have
\beq
0<\theta_1=\dfrac{\pi}{2g}+\theta<\dfrac{\pi}{g},\quad \cot\dfrac{\pi}{g}<\la_1(\theta)<\infty.
\label{0theta1}
\enq

\noindent
{\bf Fact 3.}
For any normal geodesic $\ga_\theta$ of $M_\theta$, $M_\theta\cap\ga_\theta$ is a parallel $2g$-gon, 
i.e. a $2g$-gon deformed from a regular $2g$-gon by a parallel transformation 
in the Lie geometric sense (see Example 5.2).  %

\smallskip\noindent
{\bf Fact 4.}
If $M_\theta$ has distinct $g$ principal curvatures,  
 the multiplicity $m_i$ of $\la_i$ is common when $g=1,3,6$, and 
$m_1=m_{\text{odd}}$ and $m_2=m_{\text{even}}$ when $g=2,4$. 
In \S2 of \cite{TY}, Tang-Yan show that the mean curvature $H_\theta$ of $M_\theta$ is given by
 \beq
H_\theta=\dfrac{g}2\Big(m_1t-\dfrac{m_2}{t}\Big),\quad t=\cot\dfrac{g\theta_1}2, 
\quad 0<\theta_1=\dfrac{\pi}{2g}+\theta<\dfrac{\pi}{g}. 
\label{HH4}
\enq
Thus $H_\theta$ is monotone w.r.t. $\theta$, and $\theta$ is uniquely determined by 
 the mean curvature. 
The minimal case occurs when 
\beq
t^2=\cot^2\dfrac{g\theta_1}2=\dfrac{m_2}{m_1}. 
\label{HH}
\enq
In particular when $g=1,3,6$,  
$m_1t-\dfrac{m_2}{t}=m\dfrac{t^2-1}{t}=2m\cot g\theta_1$, and so
 \beq
H_\theta=gm\cot(g\theta_1) 
\label{HHm}
\enq
holds, which vanishes only when $\theta=0$, i.e.,  $\theta_1=\dfrac{\pi}{2g}$. 
When $g=2,4$, $M_0$ is not necessarily minimal, 
 since $H_\theta$ depends on $m_1,m_2$ in (\ref{HH4}). 

\section{Closed case: $g=3$} 
\label{sclosedg3}

Hereafter, we assume $M$ be a {\it closed embedded\,} Dupin. 
Then $M$ is taut \cite{Th}, namely, the function (\ref{distance}) 
is generically a perfect Morse function. 
In the following, we use this important feature: tautness in the background. 

 When  $g=3$, the principal curvatures $\la_1>\la_2>\la_3$ 
have the same multiplicity $m\in \{1,2,4,8\}$ \cite{C2}. 
In the non-closed case, there is freedom for multiplicities (Example 2.1 (3)), 
and the following argument does not work. 

 For later use, we put $\la_1=\la$, $\la_2=\mu$ and $\la_3=\nu$, 
 and denote by $D_i$  the corresponding curvature distributions which decompose the tangent bundle  
  into $D_1\oplus D_2\oplus D_3$. 
 We denote any unit vector in $D_j$ by $e_j$. 
 The mean curvature is given by 
  \beq
H= m(\la+\mu+\nu),  \quad m\in \{1,2,4,8\} 
\label{g3H}
 \enq
where we do not take the mean. 
 
\bel
When $M$ is a closed CMC Dupin with $g=3$,  
there exists a point $p\in M$ at which all $d_{ji}=e_j(\la_i)$ vanish. 
 \label{ldji}
 \enl
 \Proof
As  $H$ is constant, we have for any $e_j\in D_j$, 
 \beq
 d_{j1}+d_{j2}+d_{j3}=0,\quad  j=1,2,3.
  \label{dj10}
 \enq
 Let $p\in M$ be a critical point of $\la=\la_1$ on $M$.  Then we have
 $$
 d_{j1}=0,\quad j=1,2,3.
 $$
 Since $d_{jj}=0$ as $M$ is Dupin, we have at $p$:   
  $$
 d_{12}+ d_{13}=0,\quad  d_{23}=0,\quad  d_{32}=0. 
 $$
  Thus if we choose a critical point  $p^1$ of $\la_2$ on $L_1(p)$, we have done. 
   \qed
  
\begin{figure}[htbp]
\centering
\includegraphics[width=0.4\textwidth]{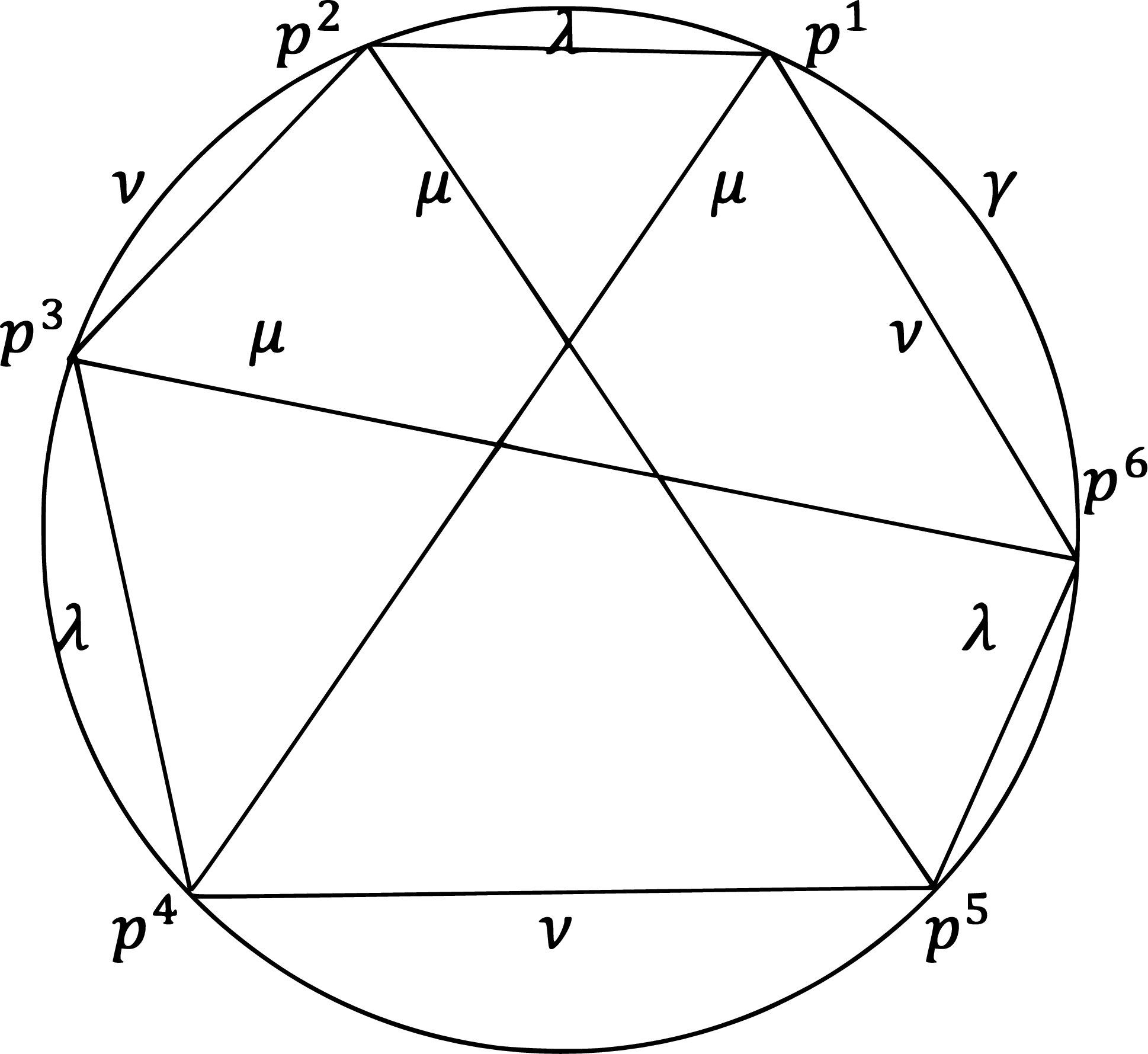}
    \caption{Link of leaves}
    \label{fg3link}
\end{figure}

\begin{figure}[htbp]
\centering
\begin{minipage}[b]{0.49\columnwidth}
    \centering
    \includegraphics[width=0.98\columnwidth]{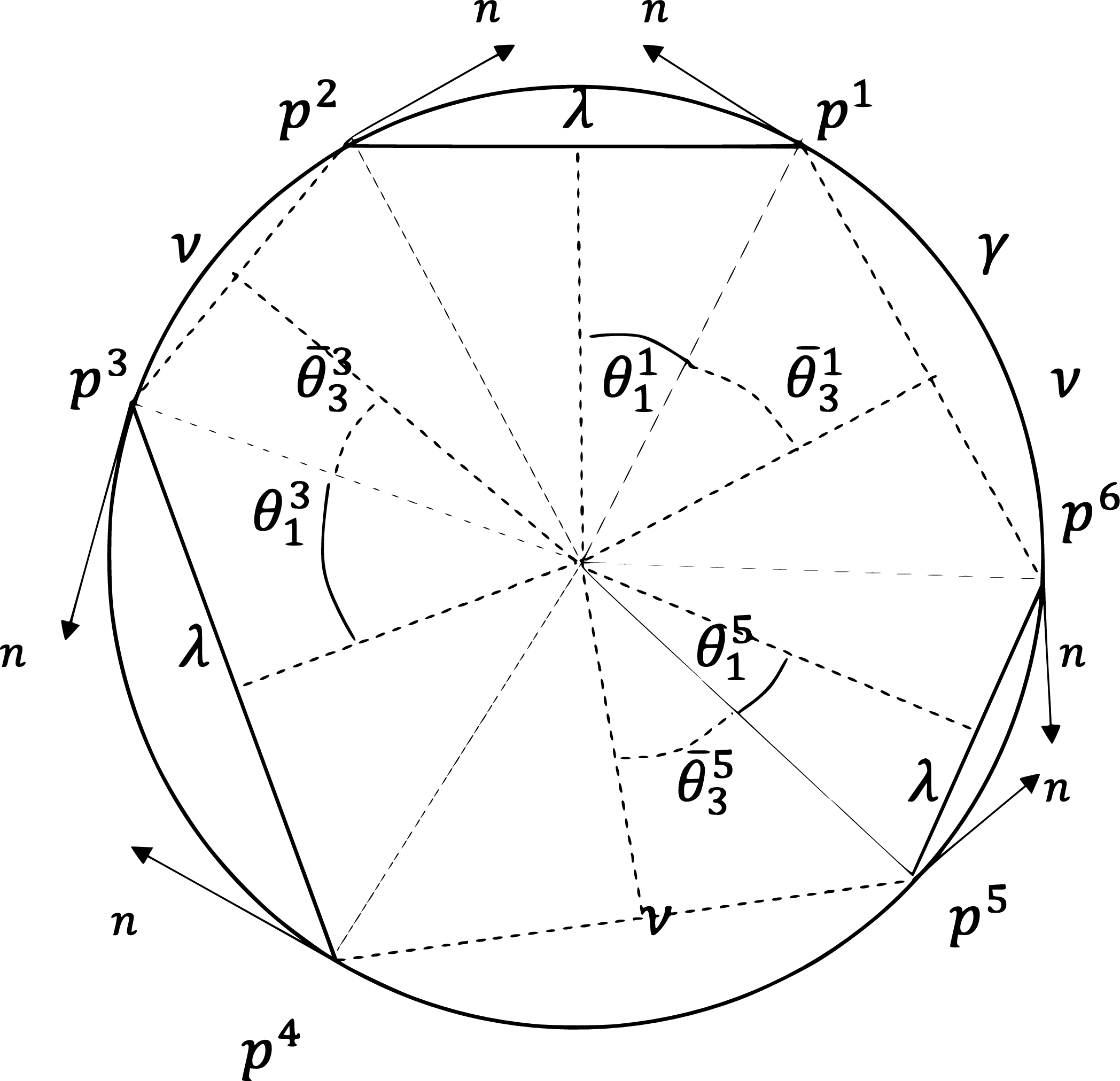}
    \caption{Angles}
    \label{fAng}
\end{minipage}
\begin{minipage}[b]{0.49\columnwidth}
    \centering
    \includegraphics[width=0.93\columnwidth]{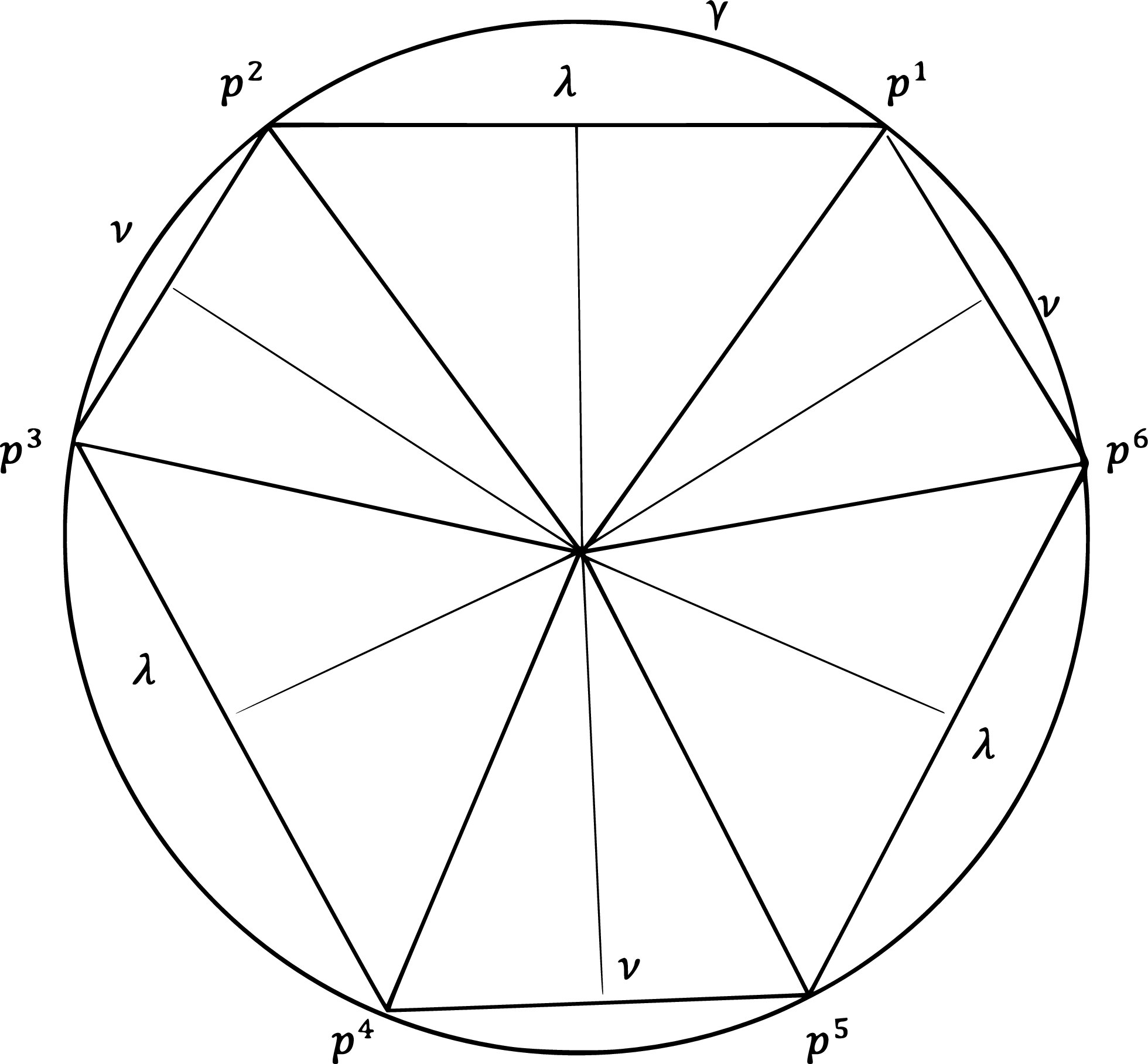}
    \caption{Parallel hexagon}
    \label{fig:b}
\end{minipage}
\end{figure}

 At $p=p^1$ in the lemma, according to the argument in \cite{M3}, each leaf $L_i(p^1)$ of $D_i$ 
 is totally geodesic in the corresponding curvature sphere $C_i$, and cuts the normal geodesic $\ga$ 
through $p^1$ orthogonally at $p^{2i}\in \ga$, $i=1,2,3$, where $p^{2i}$ is 
the antipodal point of $p^1$ in each leaf $L_i(p^1)$. 
The situation is the same as Figure 4.1 $\sim$ 4.3 in \cite{M3}, where by the argument using the tautness, 
we see $M\cap\ga=\{p^1,\dots,p^6\}$, and each leaf at $p^t$ intersects with some  other 
leaves at points on $\ga$ as in Figure \ref{fg3link}. 
More precisely, denoting $\la^t=\la(p^t)$ etc., we have 
\beq
\bary
\la^1=\la^2,\quad \la^3=\la^4,\quad \la^5=\la^6\\
\mu^1=\mu^4,\quad \mu^2=\mu^5,\quad \mu^3=\mu^6\\
\nu^1=\nu^6,\quad \nu^2=\nu^3,\quad \nu^4=\nu^5. 
\eary
\label{link}
\enq
\bel 
 Assume $\la$ takes its maximum at $p^1\in M$, and  $\la^1\ge\la^3\ge\la^5$ w.l.g. 
Then from CMC, $\mu^1\le\mu^2=\mu^5\le\mu^3$ and $\nu^3\le\nu^1\le\nu^5$ follow. 
 \label{llamunu}
 \enl
 \Proof
 Since we assume $\la^1\ge\la^3\ge\la^5$,  
 $\mu^1+\nu^1\le \mu^3+\nu^3\le \mu^5+\nu^5$ holds.  
 From $\la^1=\la^2$, we have $\mu^1+\nu^1=\mu^2+\nu^2\le \mu^3+\nu^3$, 
 and  $\nu^2=\nu^3$ implies $\mu^2\le \mu^3$. 
 Similarly from $\la^3=\la^4\ge\la^5$, 
 we have $\mu^3+\nu^3=\mu^4+\nu^4\le \mu^5+\nu^5$,  
 and  $\nu^4=\nu^5$ implies $\mu^4\le\mu^5$, i.e.,  
 $\mu^1\le\mu^2$.  
  Then   $\la^1+\nu^1\ge\la^2+\nu^2$   follows, and $\nu^1\ge\nu^2$ holds.  
  Finally from $\la^5=\la^6$,  $\mu^5+\nu^5=\mu^6+\nu^6$ follows, 
  and $\mu^5=\mu^2\le\mu^3=\mu^6$ implies $\nu^5\ge \nu^6=\nu^1$. 
   \qed
 \bep
In the situation of Lemma \ref{llamunu},  each principal curvature coincides at all $p^t$'s, 
and $M\cap\ga$ is a parallel hexagon (Figure \ref{fig:b}).  
 \label{ppconst}
 \enp
 \Proof
Put 
\beq
\la^t=\cot\theta^t_1,\quad \mu^t=\cot\theta^t_2,\quad \nu^t=\cot\theta_3^t,\quad 
0<\theta^t_i<\pi, \quad t=1,\dots,6, 
\label{tehtat}
\enq
and 
\beq
\bar\theta_3^t=\pi-\theta_3^t,
\label{bartheta3}
\enq
which satisfies $0<\bar\theta_3^t<\pi$, and  
\beq
\nu^t=\cot\theta_3^t=-\cot\bar\theta^t_3.
\label{theta33}
\enq
Recall the cotangent is a decreasing function (-cotangent is increasing).  
  
\noindent
We easily see (Figure \ref{fAng}), noting the direction of the normal, 
 \beq
 \bary
\mu^1=\cot\theta^1_2=\cot\frac12\angle(p^1Op^4)=\cot(\theta^1_1+\bar \theta^3_3+\theta^3_1),\\
\mu^2=\mu^5=\cot\theta^5_2=\cot\frac12\angle(p^5Op^2)=\cot(\theta^5_1+\theta^1_1+\bar \theta^1_3), \\
 \mu^3=\cot \theta^3_2=\cot\frac12\angle(p^3Op^6)=\cot(\theta^3_1+\bar \theta^5_3+\theta^5_1).
 \eary
   \label{thetabartheta}
   \enq
   Since $ \nu^3\le \nu^1$ ($\bar \theta^3_3\le  \bar \theta^1_3$) and $\la^3\ge\la^5$ ($\theta^3_1\le\theta^5_1$), 
   we have 
  $\bar \theta^3_3+\theta^3_1\le \bar \theta^1_3+\theta^5_1$, and $\mu^1\ge\mu^2$ follows. 
   Then by the previous lemma, we obtain $\mu^1=\mu^2$, and $\nu^1=\nu^2$ follows as $\la^1=\la^2$. 
   Also from $\la^1\ge\la^3$ and $\nu^1\le\nu^5$, 
   we have $ \theta^1_1+\bar \theta^1_3\le \theta^3_1+\bar \theta^5_3$, and $\mu^2\ge\mu^3$ follows.
   Then by the lemma, we have
     $$
   \mu^1=\mu^2=\mu^3(=\mu^6).
   $$
  Thus from $\nu^2=\nu^3$,  
   $\la^2=\la^3$ follows, and from $\nu^1=\nu^6$, $\la^1=\la^6$ follows. 
   Finally we obtain $\la^1=\la^3=\la^5$, and  $\nu^1=\nu^3=\nu^5$. 
   \qed
   
  \Remark 
$M\cap\ga$ is a parallel hexagon when $\la_1^t$ and $\la_3^t$ are independent of $t$. 

 \smallskip
 Next, let $q$ be a minimum point of $\la$, and let $q^5$ be a critical point of $\mu$ on $L_1(q)$. 
Then as before, all $d_{ji}$ vanish at $q^5$, 
and denoting the normal geodesic  through $q^5$ by $\ga'$,  
we obtain  $M\cap\ga'=\{q^1,\dots,q^6\}$ replacing $p^t$ by $q^t$ in Figure \ref{fg3link}.   
Putting $\la_q^t=\la(q^t)$ etc.,  we may assume $\la^5_q\le \la_q^3\le\la_q^1$ w.l.g. 
Hence replacing $p^t$ by $q^t$ in the above argument, we know $M\cap\ga'$ is also a 
parallel hexagon.  

\medskip\noindent
{\em Proof of Theorem \ref{tmain} (i):} 
 By Proposition \ref{ppconst}, $\{p^1,\dots, p^6\}=M\cap\ga$ is a parallel hexagon, and so 
 isometric to $M_\theta\cap\ga_\theta$ for some $\theta$. 
Here, $\theta$  is uniquely determined  
 by the mean curvature $H$ of $M$ (Fact 4). 
Thus at $p^1$, we obtain $\la=\la_\theta=\la_1(\theta)$  in (\ref{thetai}). 
Since the same is true at the minimum point $q^5$ of $\la$, $\la$ is constant all over $M$.  
 Next consider $p$ and $q$ at which $\max\mu$ and $\min\mu$ is taken, resp. 
As we have $d_{j1}=0$ and $d_{12}=0$ at $p,q$,  
 a similar argument  
 implies $M\cap\ga$ and $M\cap\ga'$ at $p$ and $q$ are isometric to $M_\theta\cap\ga_\theta$,  
and $\max\mu=\mu_\theta=\min\mu$ follows, namely, $\mu$  is constant all over $M$. 
As $M$ has CMC, $\nu$ is constant, and  $M$ is isoparametric. 
 \qed
 
 \medskip

\Remark
Tang-Yan  \cite{TY} give the scalar curvature of $M_\theta$ : 
\beq
R_\theta =9m(m-1)\big(1+\cot^2(3\theta_1)\big),\quad \theta_1=\dfrac{\pi}6+\theta. 
\label{scalar3}
\enq
Thus when $m=1$ $M_\theta$ is scalar flat, and when $m>1$, scalar positive.

\section{Closed CMC and CSC Dupin with $g=4$} 
\label{sg4}

When $g=4$, let $\la_1>\la_2>\la_3>\la_4$ be the principal curvatures 
with multiplicities $m_1=m_3,m_2=m_4$. 
Let $D_i$ be the curvature distribution for $\la_i$. 
As before, we denote any unit vector in $D_j$ by $e_j$. 
We show

\smallskip\noindent
{\bf Theorem \ref{tmain} } (ii). 
{\em A closed CMC Dupin hypersurface in $S^n$ with $g=4$ is isoparametric if it has 
 constant scalar curvature \textup{(CSC)}.}

\smallskip
Since the scalar curvature $R$ is given by
\beq
R=(n-1)(n-2)+H^2-\|A\|^2,
\label{RR}
\enq
$H$ and $R$ are constant when  
\beq
\bary
m_1\la_1+m_2\la_2+m_3\la_3+m_4\la_4=H\\
m_1\la^2_1+m_2\la^2_2+m_3\la^2_3+m_4\la^2_4=S 
\label{ms4}
\eary
\enq
are constant. 
In this section, we always assume CMC and CSC. 

\Remark
Hence two principal curvatures determine the other two. 
\label{Rem 4.1}

\smallskip
For $e_j\in D_j$, 
we have  
\beq
\sum_{i =1}^4m_id_{ji}=0,\quad \sum_{i=1^4}m_i\la_id_{ji}=0, \quad j=1,\dots 4, 
\label{dij2}
\enq
namely, say, for $j=1$, 
$$
\bary
m_2d_{12}+m_3d_{13}+m_4d_{14}=0\\
m_2\la_2d_{12}+m_3\la_3d_{13}+m_4\la_4d_{14}=0,
\eary
$$
because $e_1(\la_1)=0$  as $M$ is Dupin. 
Thus $d_{13}$ and $d_{14}$ are obtained from $d_{12}$.  
Especially, they vanish when $d_{12}=0$, since $m_3m_4\la_4-m_3m_4\la_3\ne0$. 
Just in the same way, we have 
\bel
For distinct $i,j,k,l\in\{1,2,3,4\}$, $d_{jk}$ and $d_{jl}$ are obtained from $d_{ji}$. 
Especially,  $d_{ji}=0$ implies $d_{jk}=d_{jl}=0$. 
\label{lijk}
\enl

\bel There exists a point $p\in M$ at which  all $d_{ji}$ vanish. 
\label{ldji04}
\enl
\Proof
Let $p\in M$ be a critical point of $\la_1$. Then on $L_1(p)$,  $d_{j1}=0$ holds, 
and by Lemma \ref{lijk}, $d_{j2}=d_{j3}=d_{j4}=0$ holds for $j=2,3,4$. 
Next consider $d_{1i}$, $i=2,3,4$. 
Since $\la_2$ is critical  at some point of $L_1(p)$, take such point $p^1\in L_1(p)$.  
Then $d_{12}=0$, implying $d_{13}=d_{14}=0$ by the same lemma. 
 Thus all $d_{ji}$'s vanish at $p^1$. 
\qed

Let  $\ga$ be the normal geodesic of $M$ at $p^1$. 
Since the leaf $L_i$ of $D_i$ is totally geodesic in the curvature 
sphere $C_i$ if and only if  $\la_i$ takes its critical value on $M$ (\cite{M3}, \cite{M4}), 
$L_i(p^1)$  is totally geodesic in $C_i$, and we may put  $L_1(p^1)\cap \gamma=\{p^1,p^2\}$, 
$L_2(p^1)\cap\gamma=\{p^1,p^4\}$, $L_3(p^1)\cap\gamma=\{p^1,p^6\}$ and 
$L_4(p^1)\cap\gamma=\{p^1,p^8\}$, 
where $p^{2i}$ is the antipodal point of $p^1$ in $L_i(p^1)$. 

\begin{figure}[htbp]
\centering
\begin{minipage}[b]{0.40\columnwidth}
    \centering
    \includegraphics[width=1.0\columnwidth]{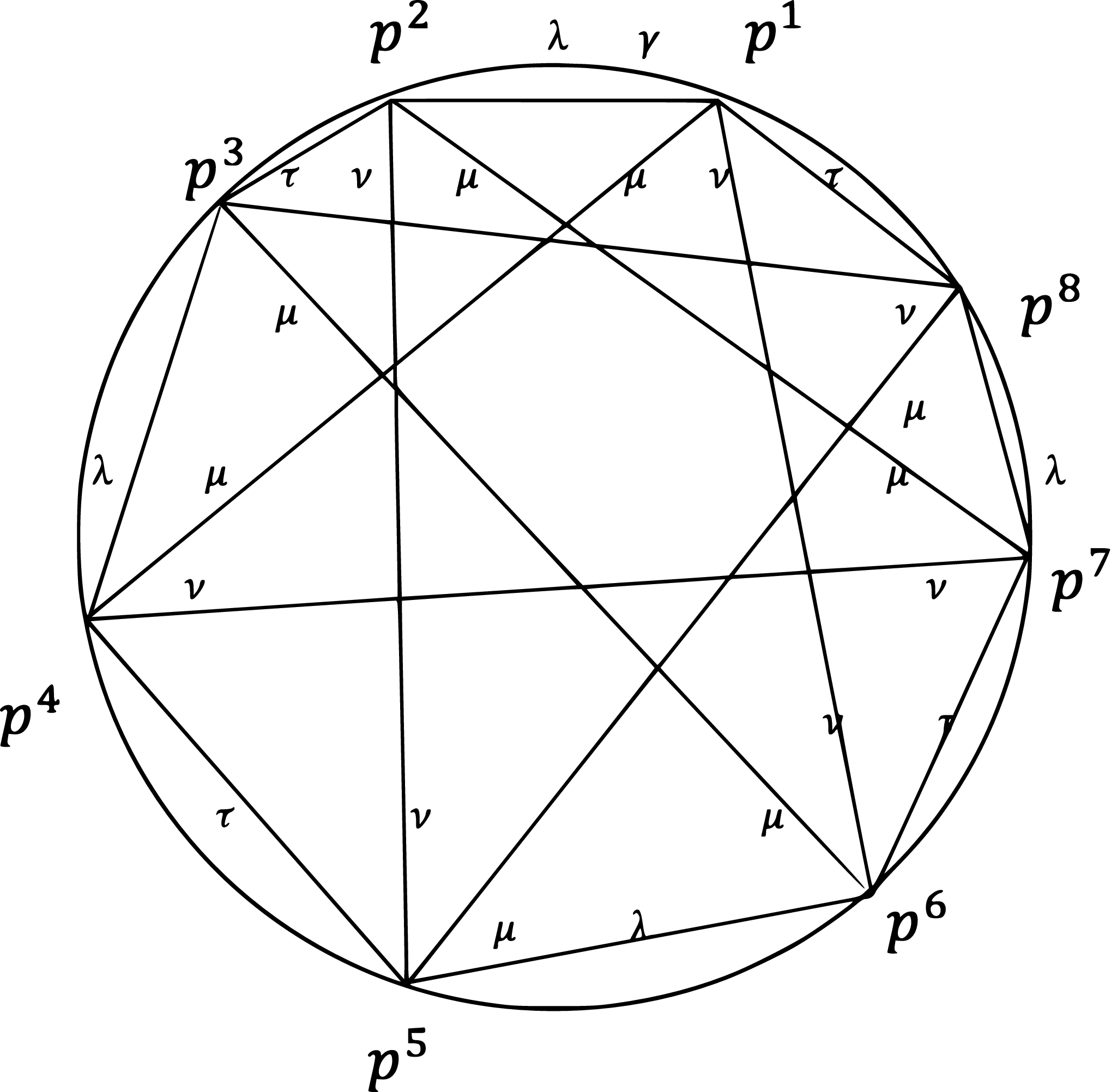}
    \caption{Link of leaves}
    \label{fg4link}
\end{minipage}
\begin{minipage}[b]{0.48\columnwidth}
    \centering
    \includegraphics[width=1.0\columnwidth]{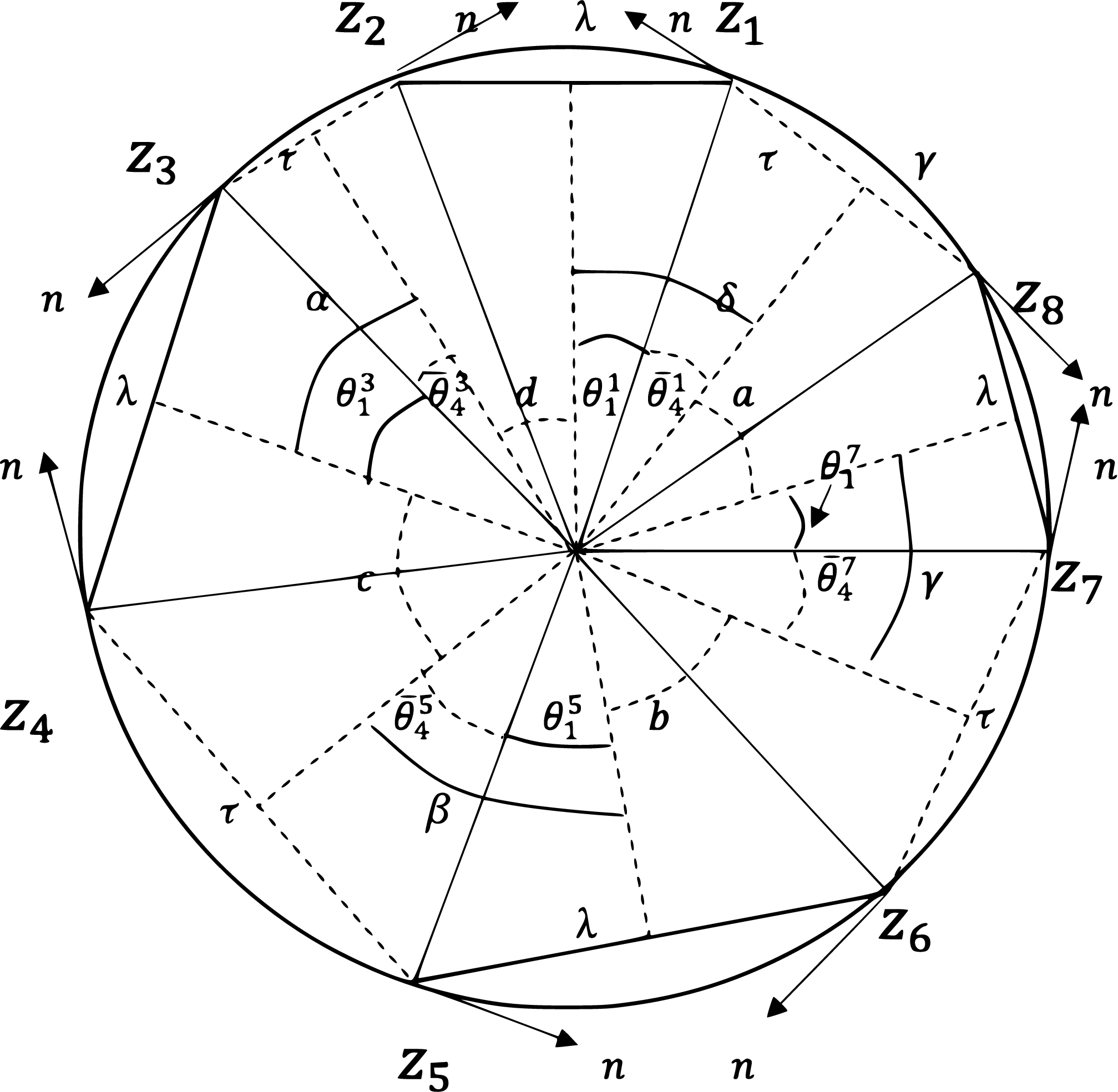}
    \caption{Angles}
    \label{fg4Ang}
\end{minipage}
\end{figure}

 Thus, at $p^1$, the situation becomes the same as in Proposition 6.1 of  \cite{M4}.   
Applying this proposition where tautness is essential, we obtain:  
\bep
We have $M\cap\gamma=\{p^1,p^2,p^3,p^4, p^5,p^6,p^7,p^8\}$ 
where $p^1,p^2,p^4,p^6,p^8$  
 are as above, and $p^3$, $p^5$, $p^7$ are as  in Figure \ref{fg4link}, 
denoting $\la=\la_1$, $\mu=\la_2$, $\nu=\la_3$ and $\tau=\la_4$. 
The leaves at each $p^t$ are all totally geodesic 
in each curvature sphere and cut $\gamma$ orthogonally at some $p^s$. 
They intersect with each other as in Figure \ref{fg4link}.    
\label{pint}
\enp
 
\bep
Denote $\la_i^t=\la_i(p^t)=\cot \theta_i^t$, $0<\theta_i^t <\pi $,  for $i=1,2,3,4$, and $t=1,\dots,8$.  
Putting 
\beq
\bary
\theta^1_2=\theta^1_1+\al,\, \,\theta^1_3=\theta^1_2+\be, \,\, 
\theta^1_4=\theta^1_3+\ga\\ 
~\\
\theta^2_1=\theta^1_1,\,\, \theta^2_2=\theta^1_1+a,\,\, \theta^2_3=\theta^2_2+b,\,\,
\theta^2_4=\theta^1_3 +c, 
\eary
\enq
$\de=\pi-(\al+\be+\ga)$ and $d=\pi-a-b-c$, we have 
the following table of the relation among angles: 

\begin{center}
\begin{tabular}{|c|r|r|r|r|r|r|r|r|c|} \hline
point & $\theta_1^t$ &  & $\theta_2^t$ & & $\theta_3^t$ & & $\theta_4^t $ & & $\pi+\theta^t_1 $\\ \hline
 $p^1$ &$\theta_1^1$  &$+$ & $\al$& $+$ &$\be$& $+$ &$\ga$&$+$ &$\de$\\
 $p^2$ &$\theta_1^1$  &$+$&$a$  & $+$ &$b$& $+$&$c$ &$+$ &$d$\\
$p^3$ &$\theta_1^3$ &$+$ &$\be$& $+$ &$\ga$&$+$  &$\de$&$+$  &$\al$\\
$p^4$ &$\theta_1^3$ &$+$ &$d$& $+$&$a$ & $+$ &$b$&$+$ &$c$ \\
$p^5$ &$\theta_1^5$ &$+$ &$\ga$& $+$ &$\de$&$+$  &$\al$& $+$&$\be$\\
$p^6$ &$\theta_1^5$ & $+$ &$c$&$+$ &$d$ & $+$ &$a$&$+$  &$b$\\
$p^7$ &$\theta_1^7$ &$+$ &$\de$& $+$ &$\al$&$+$  &$\be$&$+$ &$\ga$\\ 
$p^8$ &$\theta_1^7$ &$+$ &$b$& $+$ &$c$&$+$ & $d$&$+$ &$a$\\ \hline
 \end{tabular}
 \captionof{table}{Angles} 
 \end{center}
 \label{palbega4}
\enp

\Proof 
Consider $\ga$ as a unit circle. 
Noting the direction of the normal given in Figure \ref{fg4Ang} 
(where we denote $p^t$ by $z^t$),  
we have as in (\ref{thetabartheta}), 
$$
\bary
\mu^1=
\cot(\theta^1_1+\theta^3_1+\bar \theta^3_4),\,
\nu^1=\cot(\theta^1_2+\theta^5_1+\bar\theta^5_4),\,
\tau^1=\cot(\theta^1_3+\theta^7_1+\bar\theta^7_4)\\
\mu^2=
\cot(\theta^1_1+\theta^8_1+\bar \theta^8_4),\,
\nu^2=\cot(\theta^1_2+\theta^6_1+\bar\theta^6_4),\,
\tau^1=\cot(\theta^1_3+\theta^4_1+\bar\theta^4_4)\\
\eary$$
and so  
$$
\bary
\al=\theta^3_1+\bar\theta^3_4,\quad \be=\theta^5_1+\bar\theta^5_4,\quad 
\ga=\theta^7_1+\bar\theta^7_4\\
a=\theta^8_1+\bar\theta^8_4,\quad b=\theta^6_1+\bar\theta^6_4,\quad 
c=\theta^4_1+\bar\theta^4_4.
\eary
$$
Also from Figure \ref{fg4link} and Figure \ref{fg4Ang}, we see,  
$$
\bary
\theta^3_2=\theta^3_1+\theta^5_1+\bar\theta^5_4=\theta^3_1+\be,\,
\theta^3_3=\theta^3_2+\theta^7_1+\bar\theta^7_4=\theta^3_2+\ga,\\
\theta^3_4=\theta^3_3+\theta^1_1+\bar\theta^1_4=\theta^3_3+\de
\eary
$$
where $\de=\theta^1_1+\bar\theta^1_4=\theta^1_1+(\pi-\theta^1_4)=\pi-(\al+\be+\ga)$. 
Since $\theta^4_1=\theta^3_1$, we have in the same way,  
$$
\bary
\theta^4_2=\theta^4_1+\theta^2_1+\bar\theta^2_4=\theta^3_1+d,\, 
\theta^4_3=\theta^4_2+\theta^8_1+\bar\theta^8_4=\theta^4_2+a,\\
\theta^4_4=\theta^4_3+\theta^6_1+\bar\theta^6_4=\theta^6_3+b
\eary
$$
where $d=\theta^2_1+\bar\theta^2_4=\theta^2_1+(\pi-\theta^2_4)=\pi-a-b-c$. 
 Similarly, $\theta^5_i$, $\theta^6_i$, $\theta^7_i$ and $\theta^8_i$ are obtained. 
\qed

\smallskip
Recall (Figure \ref{fg4link}) 
\beq
\bary
\la^1=\la^2,\quad \la^3=\la^4,\quad \la^5=\la^6,\quad \la^7=\la^8\\
\mu^1=\mu^4,\quad \mu^2=\mu^7,\quad \mu^3=\mu^6,\quad \mu^5=\mu^8\\
\nu^1=\nu^6,\quad \nu^2=\nu^5,\quad \nu^3=\nu^8,\quad \nu^4=\nu^7\\
\tau^1=\tau^8,\quad \tau^2=\tau^3,\quad \tau^4=\tau^5,\quad \tau^6=\tau^7. 
\eary
\label{relation}
\enq

\bel
We have
\beq
\theta^3_1=\theta^1_1+\al-d,\quad \theta^5_1=\theta^1_1+a+b-\ga-\de,\quad \theta^7_1=\theta^1_1+a-\de.
\label{theta31}
\enq
\label{ltheta31}
\enl
\Proof
Since $\mu^1=\mu^4$, from Table 1 follows 
$\theta^1_1+\al=\theta^3_1+d$, and 
\beq
 \theta^3_1=\theta^1_1+\al-d.
 \label{31}
 \enq
 From $\mu^2=\mu^7$, $\theta^1_1+a=\theta^7_1+\de$ follows, i.e.,  
\beq
\theta^7_1=\theta^1_1+a-\de,
\label{71}
\enq
 and from $\mu^5=\mu^8$, $\theta^5_1+\ga=\theta^7_1+b=\theta^1_1+a-\de+b$ follows and 
 \beq
\theta^5_1=\theta^1_1+a+b-(\ga+\de).
\label{theta51}
\enq
\qed

\bel
Assume at $p^1$, $\la$ takes its maximum and $\mu^1\ge\mu^2$. 
Then we have at $p^1$, 
\beq
\nu^7\ge\nu^5\ge\nu^1\ge\nu^3,\quad \tau^1\ge\tau^3,\quad \tau^7\ge \tau^5.
\label{nutau}
\enq
\label{lnutau}
\enl
\Proof
Since $m_1=m_3$ and $m_2=m_4$, we have at $p^1$ and $p^8$, 
$$
\bary
\begin{cases}
m_1(\la^1+\nu^1)+m_2(\mu^1+\tau^1)=H\\ 
m_1(\la^8+\nu^8)+m_2(\mu^8+\tau^8)=H
\end{cases}\\
\begin{cases}
m_1((\la^1)^2+(\nu^1)^2)+m_2((\mu^1)^2+(\tau^1)^2)=S\\
m_1((\la^8)^2+(\nu^8)^2)+m_2((\mu^8)^2+(\tau^8)^2)=S
\end{cases}
\eary
$$
where  $\tau^1=\tau^8$. 
Taking the difference of the former two, and the latter two, respectively, we obtain
$$
\bary
m_1(\la^1-\la^8+\nu^1-\nu^8)+m_2(\mu^1-\mu^8)=0\\
m_1((\la^1)^2-(\la^8)^2+(\nu^1)^2-(\nu^8)^2)+m_2((\mu^1)^2-(\mu^8)^2)=0.
\eary
$$
Then multiplying the former  by $(\mu^1+\mu^8)$ and subtracting the second, we obtain 
\beq
(\la^1-\la^8)(\mu^1+\mu^8-\la^1-\la^8)+(\nu^1-\nu^8)(\mu^1+\mu^8-\nu^1-\nu^8)=0,  
\label{la1la8}
\enq
where $\mu^1+\mu^8-\la^1-\la^8< 0$ and $\mu^1+\mu^8-\nu^1-\nu^8>0$. 
Since $\la^1$ is max, (\ref{la1la8}) implies 
\beq 
\la^1\ge \la^8=\la^7,\quad \nu^1\ge \nu^8=\nu^3. 
\label{nu13}
\enq
A similar argument at $p^2$ and $p^7$ using $\mu^2=\mu^7$  implies 
$$
\bary
m_1(\la^2-\la^7+\nu^2-\nu^7)+m_2(\tau^2-\tau^7)=0\\
m_1((\la^2)^2-(\la^7)^2+(\nu^2)^2-(\nu^7)^2)+m_2((\tau^2)^2-(\tau^7)^2)=0, 
\eary
$$
and multiplying the former by $(\tau^2+\tau^7)$ and subtracting the second, we have
\beq
(\la^2-\la^7)(\tau^2+\tau^7-\la^2-\la^7)+(\nu^2-\nu^7)(\tau^2+\tau^7-\nu^2-\nu^7)=0,   
\label{la2la7}
\enq
where $\tau^2+\tau^7-\la^2-\la^7<0$ and $\tau^2+\tau^7-\nu^2-\nu^7<0$.  
Thus in (\ref{la2la7}), $\la^2=\la^1\ge\la^7$ implies  
\beq
\nu^2=\nu^5\le \nu^7.
\label{nu5nu7}
\enq
Next, at $p^1$ and $p^2$, using $\la^1=\la^2$, we have
$$
\bary
m_1(\nu^1-\nu^2)+m_2(\mu^1-\mu^2+\tau^1-\tau^2)=0\\
m_1((\nu^1)^2-(\nu^2)^2)+m_2((\mu^1)^2-(\mu^2)^2+(\tau^1)^2-(\tau^2)^2)=0,
\eary
$$
and multiplying the former by $(\nu^1+\nu^2)$  and subtracting the second, we have 
\beq
(\mu^1-\mu^2)(\nu^1+\nu^2-\mu^1-\mu^2)+(\tau^1-\tau^2)(\nu^1+\nu^2-\tau^1-\tau^2)=0,   
\label{mu1mu2}
\enq
where $\nu^1+\nu^2-\mu^1-\mu^2< 0$ and $\nu^1+\nu^2-\tau^1-\tau^2>0$. 
Since we have chosen $p^1$ so that $\mu^1\ge\mu^2$, we obtain from (\ref{mu1mu2}) 
\beq 
\mu^1\ge\mu^2,\quad \tau^1\ge\tau^2=\tau^3. 
\label{tau13}
\enq
This impies  $\mu^1+\tau^1\ge\mu^2+\tau^2$ and so $\la^1+\nu^1\le\la^2+\nu^2$.  
Thus   
\beq
\la^1=\la^2,\quad \nu^1\le\nu^2=\nu^5 
\label{nu12}
\enq
hold. 
Next, at $p^4$ and $p^7$, using $\nu^4=\nu^7$, we have
$$
\bary
m_1(\la^4-\la^7)+m_2(\mu^4-\mu^7+\tau^4-\tau^7)=0\\
m_1((\la^4)^2-(\la^7)^2)+m_2((\mu^4)^2-(\mu^7)^2+(\tau^4)^2-(\tau^7)^2)=0,
\eary
$$
and multiplying  the former by $(\la^4+\la^7)$ and subtracting the second, we obtain
\beq
(\mu^4-\mu^7)(\la^4+\la^7-\mu^4-\mu^7)+(\tau^4-\tau^7)(\la^4+\la^7-\tau^4-\tau^7)=0,  
\label{nu4nu7}
\enq
where $\la^4+\la^7-\mu^4-\mu^7>0$ and $\la^4+\la^7-\tau^4-\tau^7>0$. 
Since we are assuming $\mu^1=\mu^4\ge\mu^2=\mu^7$, it holds
\beq 
 \tau^4=\tau^5\le \tau^7=\tau^6,
\label{nu27}
\enq
and the lemma is proved. 
\qed

\bep
We have
\beq
\al=a,\quad \be=b, \quad \ga=c,\quad \de=d. 
\label{abcd}
\enq
\label{palabeb}
\enp
\Proof
Recall the cotangent is decreasing. By Lemma \ref{lnutau}, we have $\mu^1\ge\mu^2$ which implies
\beq
\al\le a. 
\label{ala}
\enq
The second inequality of (\ref{nutau}) is written as, using (\ref{theta31}), 
\beq
(\theta^1_1+a-\de)+\de+\al\le\theta^1_1+a+b\le\theta^1_1+\al+\be\le(\theta^1_1-d+\al)+\be+\ga
\enq
and we have immediately
\beq
\al\le b,\quad a+b\le\al+\be, \quad d\le \ga.
\label{alb}
\enq
Then from $\tau^7\ge\tau^5$,  we have
\beq
\bary
\quad(\theta^1_1+a-\de)+\de+\al+\be=\theta^1_1+a+\al+\be\\
\qquad\le(\theta^1_1+a+b-\ga-\de)+(\ga+\de+\al)
=\theta^1_1+a+b+\al, 
\eary
\enq
and so
\beq
\be\le b.
\label{beb}
\enq
Thus from (\ref{ala}) +  (\ref{beb}), and the middle innequality of (\ref{alb}),  we obtain
\beq
\al+\be=a+b,
\label{albeab}
\enq
and hence 
\beq
\al=a,\quad \be=b.
\label{abalbe}
\enq
From  $\theta_1^1=\theta_1^2$ and  (\ref{abalbe}),   $\mu^1=\mu^2$ and $\nu^1=\nu^2$ follow.  
These  imply $\tau^1=\tau^2$, and we obtain $\ga=c$ and $\de=d$.
\qed

\bel We have
\beq
\de\le \ga \le \al\le\be.
\label{ineq}
\enq
\label{lineq}
\enl
\Proof
As we assume $\la^1\ge\la^3=\cot(\theta^1_1+\al-\de)$, we have
$$
\de\le\al.
$$
From $\nu^7\ge\nu^5\ge\nu^1\ge\nu^3$ 
follows $\theta^7_3\le \theta^5_3\le \theta^1_3\le \theta^3_3$, 
and  using (\ref{theta31})  and (\ref{abcd}), we have
$$
\bary
\theta^7_1+\de+\al&=\theta^1_1+2\al\\
&\le\theta^5_1+\ga+\de=\theta^1_1+\al+\be=\theta^1_3\\
& \le\theta^3_1+\be+\ga=(\theta^1_1-\de+\al)+\be+\ga,
\eary
$$
namely, 
$$
\de\le\al\le\be,\quad \de\le \ga.
$$
Since we have  $\la^3=\la^4$,  $\nu^3\le\nu^7=\nu^4$, and $\mu^3=\cot(\theta^3_1+\be)\le \cot(\theta^3_1+\de)=\mu^4$, 
 we obtain $\tau^3\ge\tau^4$. 
 Hence $\be+\ga+\de\le \de+\al+\be$ implies  
$$
\ga\le\al. 
$$
\qed

\bep
All the principal curvatures coincide at $p^t$ for $1\le t\le 8$, and $M\cap\ga$ is a parallel octagon. 
\label{pconstpc4}
\enp
\Proof 
In addition to  $\la^3\le\la^1$, we have, using (\ref{abcd}) and (\ref{ineq}), 
$$
\mu^3=\cot(\theta^3_1+\be)=\cot(\theta^1_1+\al-\de+\be))\le \cot(\theta^1_1+\al)=\mu^1.  
$$
and  
$$
\nu^3=\cot(\theta^3_1+\be+\ga)=\cot(\theta^1_1+\al-\de+\be+\ga)\le\nu^1.  
$$
However since 
$$
\theta^3_4=\theta^3_1+ \be+\ga+\de=(\theta^1_1-\de+\al)+ \be+\ga+\de=\theta^1_1+\al+\be+\ga=\theta_4^1
$$
by the CMC condition, we must have  $\la^1=\la^3$, $\nu^1=\nu^3$ and $\mu^3=\mu^1$.  
These imply 
\beq
\al=\be=\ga=\de=\dfrac{\pi}4,
\enq
and $\theta^5_1=\theta^7_1=\theta^1_1$ follows. By Proposition \ref{palabeb} and by Table 1, 
all the principal curvatures at $p^t$ coincide for $1\le t\le 8$. 
Thus, $M\cap\ga$ is a parallel octagon. 
\qed

\begin{figure}[htbp]
\centering
\includegraphics[width=0.4\textwidth]{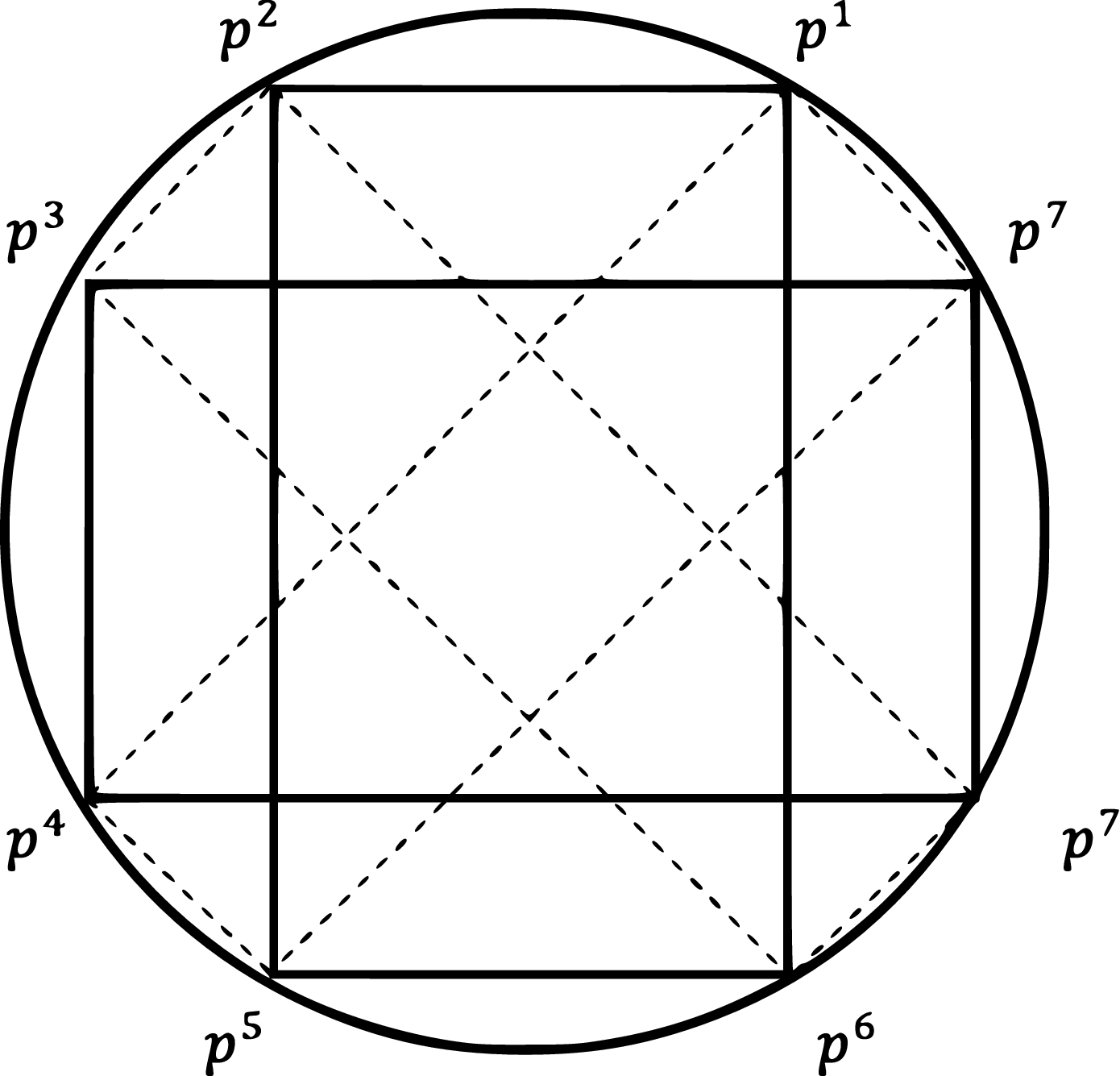}
\caption{Parallel octagon}
\label{fparag4}
\end{figure}

\smallskip
Next, we apply the above argument at a minimum point $q$ of $\la$, instead of the maximum point $p$.  
Let $q^1$ be a minimum point of $ \mu$  on $L_1(q)$. 
Then we can show all $d_{ij}$ vanish  at $q^1$ using Lemma \ref{lijk}. 
Denote by $\ga'$ the normal geodesic through $q^1$. 
Then using the tautness again, we see 
$M\cap\ga'=\{q^1,q^2,\dots,q^8\}$ and we have a similar figure as Figure \ref{fg4Ang} along $\ga'$. 

The following argument is almost parallel as before, but we need to check the change of inequalities. 
In all the argument of Lemma \ref{lnutau}, \ref{palabeb} and Lemma \ref{lineq} 
using Table 1, we replace $p^t$ by $q^t$. 
Then Lemma \ref{ltheta31} holds for $q^t$, 
and we show instead of Lemma \ref{lnutau}:
\bel
Take $q^1$ as a minimum point of $\la$ on $M$ and also minimum of $\mu$ on $L_1(q^1)$. 
Then we have, denoting $\nu^t=\nu(q^t)$, etc., 
\beq
 \nu^7\le\nu^5\le\nu^1\le\nu^3,\quad \tau^1\le\tau^3,\quad \tau^7\le \tau^5.
\label{nutauq}
\enq
\label{lnutauq}
\enl
 \Proof
As before, (\ref{la1la8}) holds at $q^1$ and $q^8$, 
 where $\la^1-\la^8$ and $\nu^1-\nu^8$ should have the same sign. 
 Since $\la^1$ is the minimum, we have 
 \beq 
 \la^1\le \la^8=\la^7,\quad \nu^1\le \nu^8=\nu^3. 
\label{nu13q}
\enq
 At $q^2$ and $q^7$ we have  (\ref{la2la7}), 
and $\la^2=\la^1\le\la^7$ implies  
\beq
\nu^2=\nu^5\ge \nu^7.
\label{nu5nu7q}
\enq
Next, at $q^1$ and $q^2$ in (\ref{mu1mu2}),  
since $q^1$ is chosen so that $\mu^1\le\mu^2$, we obtain  
\beq 
\mu^1\le\mu^2,\quad \tau^1\le\tau^2=\tau^3. 
\label{tau13q}
\enq
This impies  $\mu^1+\tau^1\le\mu^2+\tau^2$ and so $\la^1+\nu^1\ge\la^2+\nu^2$.  
As  $\la^1=\la^2$,  
\beq
 \nu^1\ge\nu^2=\nu^5 
\label{nu12q}
\enq
follows. 
Next, at $q^4$ and $q^7$ in (\ref{nu4nu7}),
from  $\mu^1=\mu^4\le\mu^2=\mu^7$, it follows
\beq 
 \tau^4=\tau^5\ge \tau^7=\tau^6.
\label{nu27q}
\enq
\qed

\bep
Under the  situation in Lemma \ref{lnutauq}, denoting the angles with $'$ in this case, we have 
\beq
\al'=a',\quad \be'=b', \quad \ga'=c',\quad \de'=d'. 
\label{abcdq1}
\enq
of which values are all $\dfrac{\pi}4$.
\label{palabebq}
\enp
\Proof
The proof is parallel to that of Proposition  \ref{palabeb}, Lemma \ref{lineq} and Proposition \ref{pconstpc4}, where 
we reverse all inequality signs and apply a similar argument. 
We check all the processes but omit details. 
\qed

Finally, we obtain:
\bep
All the principal curvatures coincide at $q^t$ for $1\le t\le 8$, and $M\cap\ga'$ is a parallel octagon. 
\label{pconstpc4q}
\enp

\medskip
{\em Proof of Theorem \ref{tmain} (ii):} 
Both octagons $p^1\dots p^8$ and  $q^1\dots q^8$ are  isometric to a parallel octagon $M_\theta \cap\ga_\theta$, 
where $\theta$ is uniquely determined by $H$ (Fact 4). 
This implies $\la(p^1)=\la_\theta=\la(q^1)$, namely, $\la$ is constant all over $M$.  
Therefore $d_{j1}=0$ holds everywhere on $M$, and 
 by Lemma \ref{lijk},  we obtain $d_{ji}=0$ for  $j=2,3,4$ and $i\ne1$. 
Taking the maximum and minimum points of $\mu$ where $d_{12}=0$ and  so 
$d_{13}=d_{14}=0$ hold, we see all $d_{ji}$ vanish at such points. 
 Then by a similar argument, we obtain $\max\mu=\min\mu$ and so $\mu$ is constant all over $M$.   
Then by Remark 4.1, $\nu$ and $\tau$ become constant.  
This proves Theorem \ref{tmain}  (ii). 
\qed

\medskip
\Remark
When $g=4$, the scalar curvature $R_\theta$ of $M_\theta$ is given by (\cite{TY})   
\beq
\bary
R_\theta &=(n-1)(n-2)+H_\theta^2-\|A_\theta\|^2\\
&=4\Big(m_1(m_1-1)(1+t^2)+m_2(m_2-1)\Big(1+\dfrac1{t^2}\Big)\Big),
\eary
\label{scalar}
\enq
where $t=\cot 2\theta_1=\cot 2\Big(\dfrac{\pi}8+\theta\Big)$,  
and  $n-1=2(m_1+m_2)$. 
Thus $M_\theta$ is scalar flat only when $m_1=m_2=1$ and otherwise, scalar positive. 
When $M_\theta$ is minimal, i.e. $t^2=\dfrac{m_2}{m_1}$ by (\ref{HH}), 
the scalar curvature is given by 
\beq
R=4(m_1+m_2)(m_1+m_2-2). 
\label{scmin}
\enq
Using the classification of isoparametric hypersurfaces with $g=4$ \cite{CCJ},\cite{Chi}, we have:  
\bep
 If  a closed {\em minimal} Dupin hypersurface $M$ in $S^{n}$ with $g=4$ has constant scalar curvature $R$, 
then $M$ is isoparametric, and $R$ lies in the discrete set 
 \beq
 \{4(m_1+m_2)(m_1+m_2-2)\}, 
 \enq
 where $(m_1,m_2)=(1,1),(2,2),(4,5)\dots$   are  infinite series given  in \textup{\cite{FKM}}.  
 The cases other than  $(m_1,m_2)=(2,2),(4,5)$  correspond to the representation of Clifford algebras. 
 \label{pscalar}
 \enp 
 
\section{Review of the Lie sphere geometry}
\label{sLiec}

Up to here, we only use elementary arguments. 
However, to prove (iii) and (iv) of Theorem\ref{tmain}, we need the Lie sphere geometry, especially the Lie curvature. 
Here, we give a brief introduction to the Lie sphere geometry. For details, see \cite{Ce}. 

Let $\R^{n+3}_2=\R^{n+1}\oplus\R^2_2$ be the space equipped with the bilinear form $\langle\,,\,\rangle_2$ 
with signature $(+\dots,+,-,-)$. 
The hyperquadrique of $\R P^{n+2}_2$ given by 
\beq
Q^{n+1}=\{[z]={}^t[z_0,\dots,z_{n+1},z_{n+2}]\in \R P^{n+2}_2\mid \langle z,z\rangle_2=0\} 
\label{QQ}
\enq
is the space of oriented hyperspheres of $S^n$. 
In fact, an oriented hyperspheres of $S^n$ centered at $p\in  S^n$ is given 
by ${}^t(p,\cos\theta,\pm\sin\theta)\in\R^{n+3}_2$, 
where $\theta$ is the radius,  and the orientation is determined by $\pm\sin\theta$, $0\le\theta\le\pi$. 
The point ${}^t(p,1,0)\in \R^{n+3}_2$ represents the point sphere, and ${}^t(p,0,1)$, 
an oriented totally geodesic hypersphere centered at $p$. 
Two elements  $k_1,k_2\in Q$ have an oriented contact if and only if  
 \beq
 \langle k_1,k_2\rangle_2=0.
 \label{contact}
 \enq
For instance, when $k_1={}^t(p,1,0)$ and $k_2={}^t(n,0,1)$ where $p$ is orthogonal to $n$, 
they satisfy (\ref{contact}), and $k_2$ has an oriented contact with $k_1$ at $p$. 
In this case,  we can define 
 a line $l=ak_1+bk_2$ of $Q$ as a one-parameter family of oriented hyperspheres having oriented contact at $p$. 
The space of lines $\Lambda^{2n-1}$ is identified with the unit tangent space $T_1S^n$, 
because $l$ is uniquely determined by a contact point $p\in S^n$ and 
the oriented unit normal $n\in T_pS^n$ to the 
contact hyperspheres at $p$. 

The Lie contact transformation group $O(n+1,2)$ is the linear group preserving 
$\langle\,,\,\rangle_2$, and  $L$  is characterized by  
 \beq
 {}^tL\bar I L=\bar I=\begin{pmatrix} I_{n+1}&0\\0&-I_2\end{pmatrix}, 
 \label{AO}
 \enq
 where  $I_{n+1}$  and $I_2$ are the unit matrices.  
 The set of column vectors (row vectors, resp.) is called a Lie frame. 
 The transformation $L$ preserves the space $Q$,  oriented contact of two elements in $Q$, 
 and  the space of lines $\Lambda^{2n-1}$. 
 
 An oriented  hypersurface $p:M\to S^n$ with unit normal vector field $n$ 
 is expressed as a {\it Lie geometric hypersurface} by the pair $(k_1,k_2)$: 
 \beq
 k_1=\bpm p\\1\\0\epm\in Q,\quad k_2=\bpm n\\0\\1\epm\in Q, 
 \label{k1k2}
 \enq
  not both $dk_1(X)$ and $dk_2(X)$ vanish for non-zero vector $X\in TM$. 
 Since $\langle k_1,k_2\rangle_2=0$, the pair $(k_1,k_2)$ defines a line $l\in  \La^{2n-1}\cong T_1 S^n$. 
 Actually, 
 \beq
 (k_1,k_2):M\to \La^{2n-1}\cong T_1 S^n
 \label{Leg}
 \enq
 is nothing but the Legendre map of $M$ into the contact manifold $T_1 S^n$. 
  
 The curvature sphere of $M$ at $p$ is an oriented hypersphere of $S^n$ contact $M$ at $p$ 
 with contact order $\ge2$. 
 For the principal curvature $\la$, the curvature sphere is given \cite{M4} by 
 \beq
v k_1+ u k_2 \in Q,\quad \la=\dfrac{v}{u}=\cot\xi,
 \label{pc}
 \enq
 where $\xi$ is the oriented radius of the curvature sphere. 
 Applying $L\in O(n+1,2)$ to $k_1,k_2$, we obtain a new Lie geometric hypersurface 
 $(Lk_1,Lk_2):M\to \La^{2n-1}$:  
 \beq
 Lk_1=\bpm q\\a\\ b \epm, \quad Lk_2=\bpm m\\c\\d\epm.
 \label{Lk}
 \enq
  The corresponding principal curvature $\tilde\la$  of the image hypersurface  is given by (4.2) in \cite{M4} as 
 \beq
\tilde\la=\dfrac{a\la +c}{b\la +d}.
 \label{newpc}
 \enq
  \Remark 
 We apply a suitable projection $\pi:\La^{2n-1}\to S^n$ to obtain a hypersurface in $S^n$ from $(Lk_1,Lk_2)$. 
Note that $q$ and $m$ themselves are not  the point  nor the normal vector of the hypersurface. 

 \Example
 A parallel transformation: 
\beq
L=\bpm L_1&0\\0&L_\theta\epm\in O(n+1)\oplus O(2), \quad 
L_\theta=\bpm \cos\theta&-\sin\theta\\\sin\theta&\cos\theta\epm, 
\label{Lpar}
\enq
 deforms a hypersurface $M$ into its parallel hypersurface $M_\theta$ (see Remark 3.8 of  \cite{M4}). 
 In fact, when $\la=\cot\xi$ is a principal curvature of $M$, 
 $M_\theta$ has the principal curvature via (\ref{newpc}), 
 \beq
  \la_\theta=\dfrac{a\la +c}{b\la +d}
  =\dfrac{ \cos\theta\cot\xi-\sin\theta}{ \sin\theta\cot\xi+\cos\theta}=\cot(\xi+\theta).
 \label{pctheta}
 \enq
 \label{Ex5.2}
 
 \smallskip
 Next, for $w_1,w_2$, $w_3,w_4\in \C$, define the cross ratio by
$$
[w_1,w_2;w_3,w_4]=\dfrac{(w_1-w_3)(w_2-w_4)}{(w_1-w_4)(w_2-w_3)}\in \C, 
$$
which takes value in $\R$ if and only if $w_1,w_2,w_3,w_4$ are concircular. 

\smallskip
 Let $\ga$ be the normal geodesic of a hypersurface $M$ in $S^n$ at $p\in M$.  
If $M$ has four principal curvatures $\la_i=\cot\theta_i$ corresponding to 
the curvature sphere $C_i$ of radius $\theta_i$,   
 $C_i$ cuts $\ga=S^1$ orthogonally at the point in the spherical distance $2\theta_i$ from $p\in M$. 
Expressing such points by $z_{2i}\in \ga=S^1\subset\C$ for later use,  
we obtain: 
\bep[Lemma 6.8 \cite{M4}]
 The Lie curvature 
 defined by  $\Phi(p)=\dfrac{(\la_1-\la_3)(\la_2-\la_4)}{(\la_1-\la_4)(\la_2-\la_3)}$ satisfies 
\beq
\Phi(p)=[z_2,z_4;z_6,z_8]\in \R, 
\label{LPhi}
\enq 
and is invariant under Lie contact transformations because of (\ref{newpc}). 
\label{pLieC}
\enp
This gives a judgment index to see whether two hypersurfaces are Lie equivalent or not. 

\section{Closed CMC and CLC Dupin with $g=4$ (iii)} 
\label{sg4iii}

In this section, we prove:

\noindent
{\bf Theorem 1.1 (iii)} 
{\em When $M$ is a closed CMC Dupin hypersurface with $g=4$, $M$ is isoparametric 
if $M$ has constant Lie curvature (CLC).}

\smallskip
The following is important in the proof: 

\smallskip\noindent
{\bf Fact 5. (Proposition 8.1 and Corollary 8.3 in \cite{M4}) } 
{\it When the Lie curvature of a Dupin hypersurface $M$ with $g=4$ is constant,  
there exists a Lie transformation at each point of $p\in M$ which maps 
$M\cap\ga$ into a regular octagon, where $\ga$ is the normal geodesic of $M$ at $p$.  }
\label{fLietr}

\smallskip
This is a local fact, and no global result comes out from this. 
However, since Lie curvatures are Lie invariant, 
the values are computed from the principal curvatures of the minimal isoparametric 
hypersurface $\bar M$ (see Fact 4), as $\bar M\cap\ga$ is a regular octagon: 
\beq
\la=\sqrt2+1=-\tau,\quad \mu=\sqrt2-1=-\nu. 
\label{pcmini46}
\enq
Even though, we cannot directly connect this fact with the mean curvature of $M$, 
because the Lie contact transformation is not compatible with the metric structure.

\subsection{Critical point of all the principal curvatures}
\label{sscritical}

\bel When $M$ is a closed CMC Dupin hypersurface with $g=4$ having CLC, 
there exists $p\in M$ at which all $d_{ji}$ vanish. 
\label{ldji0}
\enl
\Proof
Because
\beq
H=m_1\la+m_2\mu+m_1\nu+m_2\tau, \quad
\Phi=\dfrac{(\la-\mu)(\nu-\tau)}{(\la-\tau)(\nu-\mu)}=-1
\label{HPhi}
\enq
are constant on $M$ where the value of $\Phi$ is computed from (\ref{pcmini46}),  
we can describe $\mu,\tau$ by $\la,\nu$. 
Actually, putting 
\beq
A=\mu+\tau=\dfrac1{m_2}(H-m_1(\la+\nu)), 
\label{AA}
\enq
we have from the second equation of (\ref{HPhi}) 
$$
0=(\la-\mu)(\nu-\tau)+(\la-\tau)(\nu-\mu)=2\mu\tau- A(\la+\nu)+2\la\nu,
$$
and so putting 
\beq
B:=\mu\tau=\dfrac12\{A(\la+\nu)-2\la\nu\}=\dfrac12\{\dfrac1{m_2}(H-m_1(\la+\nu))(\la+\nu)-2\la\nu),
\label{BB}
\enq 
we know $\mu, \tau$ are two solutions of 
\beq
t^2-At+B=0.  
\label{AB}
\enq
Then on a $\la$-leaf $L_1$, $\mu, \tau$  are functions of only $\nu$, 
which implies 
\beq
d_{12}=f(\nu)d_{13},\quad d_{14}=g(\nu)d_{13}.
\label{d1214}
\enq
Thus if $\la$ is critical on $M$ and $\nu$ is critical at $p$ on $L_1(p)$, we have at $p$,  
\beq
d_{j1}=0,\quad d_{1i}=0,\quad 1\le j,i\le 4.
\label{dj1}
\enq
On the other hand, from (\ref{HPhi}), we obtain
$$
\bary
0&=m_1d_{j1}+m_2d_{j2}+m_1d_{j3}+m_2d_{j4}=0\\
0&=e_j(\log\Phi)=\dfrac{d_{j1}-d_{j2}}{\la-\mu}-\dfrac{d_{j1}-d_{j4}}{\la-\tau}+\dfrac{d_{j3}-d_{j4}}{\nu-\tau}
-\dfrac{d_{j3}-d_{j2}}{\nu-\mu}\\
&=\Big(\dfrac1{\nu-\mu}-\dfrac1{\la-\mu}\Big)d_{j2}+\Big(\dfrac1{\nu-\tau}-\dfrac1{\nu-\mu}\Big)d_{j3}
+\Big(\dfrac1{\la-\tau}-\dfrac1{\nu-\tau}\Big)d_{j4},
\eary
$$
and at $p$, the following hold:
\beq
\begin{cases}
m_1d_{23}+m_2d_{24}=0\\
d_{32}+d_{34}=0\\
m_2d_{42}+m_1d_{43}=0,
\end{cases}\\
\label{d230}
\enq
\beq
\begin{cases}
\dfrac{\tau-\mu}{\nu-\mu}d_{23}+\dfrac{\nu-\la}{\la-\tau}d_{24}=0\\
\dfrac{1}{(\nu-\mu)(\la-\mu)}d_{32}-\dfrac{1}{(\la-\tau)(\nu-\tau)}d_{34}=0\\
\dfrac{\la-\nu}{\la-\mu}{d_{42}}+\dfrac{\tau-\mu}{\nu-\tau}d_{43}=0. 
\end{cases}
\label{d232}
\enq
 In (\ref{d230}), the coefficients are all positive, 
and in (\ref{d232}), two coefficients of the first (the third, resp.) equation have opposite signs.  
Hence $d_{23}=d_{24}=d_{42}=d_{43}=0$ follows. 
Next, the second of (\ref{d232}) is rewritten as
$$
\bary
d_{32}- \dfrac{(\la-\mu)(\nu-\tau)}{(\la-\tau)(\nu-\mu)}\Big(\dfrac{\nu-\mu}{\nu-\tau}\Big)^2d_{34}\\
\qquad =d_{32}-\Phi\Big(\dfrac{\nu-\mu}{\nu-\tau}\Big)^2d_{34}
=(1-\Big(\dfrac{\nu-\mu}{\nu-\tau}\Big)^2\Big)d_{32}=0,
\eary
$$
where we use $\Phi=-1$ and the second equation of (\ref{d230}).
Thus from $\dfrac{\nu-\mu}{\nu-\tau}<1$, it follows $d_{32}=d_{34}=0$, 
namely, all $d_{ji}$'s vanish at $p$. 
\qed

Finally at $p=p^1$, the situation becomes the same as in Proposition 6.1 of \cite{M4}, 
stated as Proposition \ref{pint} in \S\ref{sg4} in the present paper. 

\smallskip
Next, we denote $p^1$, $p^2$ by $z_1$, $z_2$ where 
$$
z_2=e^{2i\theta_1}z_1, 
$$
and $\la^1=\cot\theta_1$. Moreover, we put
$$
z_4=e^{2i(\theta_1+\al)},\quad z_6=e^{2i(\theta_1+\al+\be)},\quad z_4=e^{2i(\theta_1+\al+\be+\ga)},
$$
where 
$$
\al=\theta^3_1+\bar\theta^3_4,\quad \be=\theta^5_1+\bar\theta^5_4,\quad 
\ga=\theta^7_1+\bar\theta^7_4,
$$
which follows from Table 1.  
Now from (\ref{HPhi}), we have
$$
\bary
-1&=\Phi=\dfrac{(\la-\mu)(\nu-\tau)}{(\la-\tau)(\nu-\mu)}
=[z_2,z_6:z_4,z_8]
=\dfrac{(z_2-z_4)(z_6-z_8)}{(z_2-z_8)(z_6-z_4)}\\
&=\dfrac{(e^{2i\theta_1}-e^{2i(\theta_1+\al)})(e^{2i(\theta_1+\al+\be)}-e^{2i(\theta_1+\al+\be+\ga)})} 
{(e^{2i\theta_1}-e^{2i(\theta_1+\al+\be+\ga)})(e^{2i(\theta_1+\al+\be)}-e^{2i(\theta_1+\al)})}\\
&=\dfrac{(1-e^{2i\al})(e^{2i(\al+\be)}-e^{2i(\al+\be+\ga)})}
{(1-e^{2i(\al+\be+\ga)})(e^{2i(\al+\be)}-e^{2i\al})}
=\dfrac{(1-e^{2i\al})(e^{2i\be}-e^{2i(\be+\ga)})}
{(1-e^{2i(\al+\be+\ga)}(e^{2i\be}-1)}, 
\eary
$$
and so 
$$
\bary
0&=(1-e^{2i\al})(e^{2i\be}-e^{2i(\be+\ga)})+(1-e^{2i(\al+\be+\ga)})(e^{2i\be}-1)\\
&=2(e^{2i\be}+e^{2i(\al+\be+\ga)})-e^{2i(\al+\be)}-e^{2i(\be+\ga)}-e^{2i(\al+2\be+\ga)}-1.
\eary
$$
Multiplying $e^{-2i\be}$ and using $\al+\be+\ga=\pi-\de$, we obtain: 
\bel
In our situation, $\Phi=-1$ implies 
\beq
2(1+e^{2i(\al+\ga)})-e^{2i\al}-e^{2i\ga}-e^{-2i\de}-e^{-2i\be}=0.
\label{alga}
\enq
\label{leal}
\enl

\subsection{Conformal transformation}
\label{ssconf4}

\begin{figure}[htbp]
\centering
\includegraphics[width=60mm]{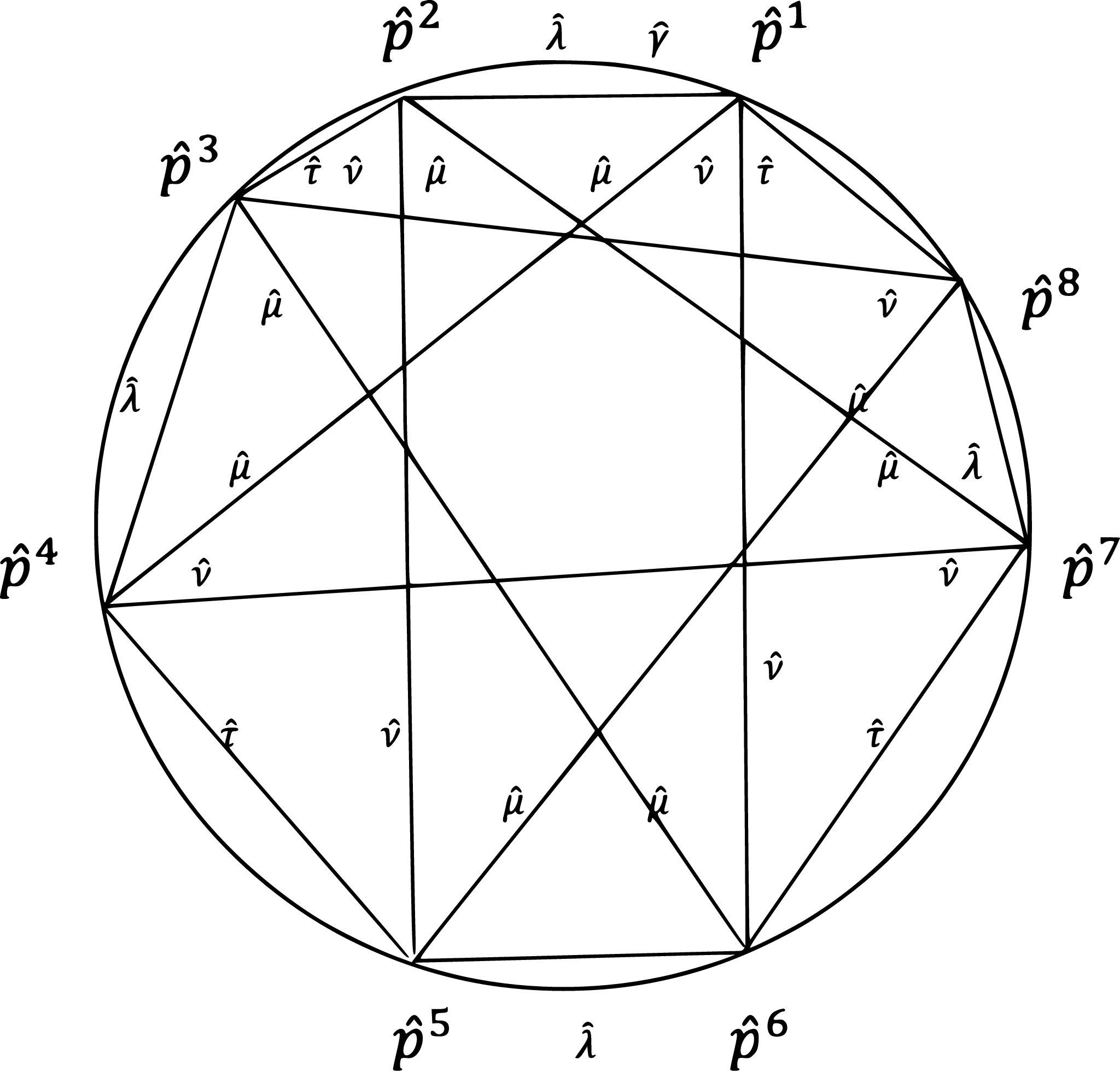}
    \caption{Conformal image}
    \label{fg4lanu}
\end{figure}

It is easy to see that there exists a conformal transformation $\hat C: M \to \hat M$, 
so that the leaves $L_\la(p^1)$ and $L_\la(p^5)$ become antipodally symmetric, 
and $L_\nu(p^1)$ and $L_\nu(p^2)$ become parallel (Figure \ref{fg4lanu}). 

We denote the image objects with hat such as $ \hat M$, $ \hat p^t$, $ \hat\la^t$ etc. 
Then the leaf $L_{\hat\la}(\hat p^1)=L_{\hat\la}(\hat p^2)$ is orthogonal to 
$L_{\hat\nu}(\hat p^1)$ at $\hat p^1$, and to $L_{\hat\nu}(\hat p^2)$ at $\hat p^2$. 
Therefore, we have, by Table 1 with hat, 
\beq
\hat\al+\hat\be=\dfrac{\pi}2=\hat a+\hat b=\hat\ga +\hat\de=\hat c+\hat d,  
\label{albepi2}
\enq
which implies, denoting 
$\hat\la^t=\cot\hat\theta^t_1$, 
\beq
\bary
 \hat\la^1\hat\nu^1=\cot\hat\theta_1\cot(\hat\theta_1+\hat\al+\hat\be)=-1,\\
 \hat\la^2\hat\nu^2=\cot\hat\theta^2_1\cot(\hat\theta^2_1+\hat a+\hat b)=-1.
\eary
\label{lanu-1}
\enq

Now apply the argument in the last subsection to $\hat M$, 
which is possible because the Lie curvature is invariant under $C$.
Then we obtain (\ref{alga}) replacing the data with that of $\hat M$, namely,  
\beq
\bary
2(1+e^{2i(\hat\al+\hat\ga)})-e^{2i\hat\al}-e^{2i\hat\ga}-e^{-2i\hat\de}-e^{-2i\hat\be}=0,\\
2(1+e^{2i(\hat a+\hat c)})-e^{2i\hat a}-e^{2i\hat c}-e^{-2i\hat d}-e^{-2i\hat b}=0
\eary
\label{alga4}
\enq
From (\ref{albepi2}),  $e^{-2i\hat\be}=e^{-2i(\pi/2-\hat\al)}=-e^{2i\hat\al}$ etc. holds,  
and (\ref{alga4}) becomes $1+e^{2i(\hat\al+\hat\ga)}=0$ and $1+e^{2i(\hat a+\hat c)}=0$.
 Thus we obtain 
\beq
\hat\al=\hat\be=\hat\ga=\hat\de=\hat a=\hat b=\hat c=\hat d=\dfrac{\pi}4.  
\label{pi4}
\enq
This means $\hat M\cap\hat \ga$ is a parallel octagon (see Figure \ref{fparag4} with hat).  
We can denote 
$$
\bary
\hat p^1=\bpm \hat u\\ \hat v\epm=-\hat p^5,\quad \hat n^1=\bpm -\hat v\\ \hat u\epm=-\hat n^5\\
\hat p^2=\bpm -\hat u\\ \hat v\epm=-\hat p^6,\quad \hat n^2=\bpm \hat v\\ \hat u\epm=-\hat n^6, \\
\hat p^3=\bpm -\hat v\\ \hat u\epm=-\hat p^7,\quad \hat n^3=\bpm -\hat u\\ -\hat v\epm=-\hat n^7\\
\hat p^4=\bpm -\hat v\\ -\hat u\epm=-\hat p^8,\quad \hat n^4=\bpm -\hat u\\ \hat v\epm=-\hat n^8.
\eary
$$
where $\hat n^t$ is the oriented unit normal of $\hat M$ at $\hat p^t$, and 
\beq
\bary
\hat\la^t=\dfrac{\hat v}{\hat u},\quad \hat\nu^t =-\dfrac{\hat u}{\hat v}, \quad \hat u^2+\hat v^2=1. \\
\eary
\label{hatlahatnu}
\enq 
More precisely, omitting $t$:  
\bel
We have $\hat\la=\cot\hat\theta_1>1$ and 
\beq
\hat\la=\dfrac{\hat v}{\hat u},\quad \hat\mu=\dfrac{\hat v-\hat u}{\hat u+\hat v},\quad
\hat\nu=-\dfrac{\hat u}{\hat v},\quad \hat\tau=\dfrac{\hat u+\hat v}{\hat u-\hat v}.
\label{hatla}
\enq
\label{lhatla}
\enl
\Proof
From (\ref{0theta1}), we have 
$\hat\theta_1=\dfrac{\pi}8+\theta<\dfrac{\pi}4$, and $\hat\la>1$ follows. 
We compute $\hat \mu$, where $\hat\tau=-\dfrac1{\hat\mu}$:
$$
\hat\mu=\cot\Big(\hat\theta_1+\frac{\pi}4\Big)=\dfrac{-1+\hat\la}{\hat\la+1}
=\dfrac{-1+\frac{\hat v}{\hat u}}{\frac{\hat v}{\hat u}+1}
=\dfrac{-\hat u+\hat v}{\hat v+\hat u}.
$$
\qed

Next, we consider the inverse $C: \hat M \to M$ of $\hat C$. 
Denote the conformal transformation $C$ restricted to $\R^4_2=\R^2\oplus\R^2_2$, where 
$\R^2$ is the 2-plane on which $\ga$ lies, by 
$$
C=\bpm *&*&*&0\\
 *&*&*&0\\
 x&y& \check\al&0\\
 0&0&0&1\epm\in O(2,1)\subset O(2,2).
 $$
Apply $C$ to $k_1=\bpm \hat p^t\\1\\0\epm$ and $k_2=\bpm \hat n^t\\0\\1\epm$, and express  
$$
Ck_1=\bpm q^t\\a_t\\0\epm, \quad Ck_2=\bpm  n^t\\c_t\\1\epm. 
$$
The principal curvature $\la_i^t$ at the original point $p^t\in M$ is given by (\ref{newpc}):   
\beq
\la_i^t= a_t\hat\la_i +c_t , \quad i=1,\dots,4, \quad t=1,\dots 8, 
\label{hatla4}
\enq
as $\hat b_t=0$ and $\hat d_t=1$. 
We know   
\beq
\bary
\begin{cases}
   a_1= \hat  ux+ \hat  vy+ \hat \al\\
   c_1=- \hat  vx+ \hat  u y\end{cases} &\quad 
\begin{cases}
   a_2=- \hat  ux+ \hat  vy+ \hat \al\\
   c_2= \hat  vx+ \hat  u y \end{cases} \\ 
\begin{cases}
a_3= -\hat  vx+ \hat  uy+ \hat \al\\
c_3= -\hat  ux-\hat v y.  \end{cases} & 
\eary
\label{hatat}
\enq
Note that the mean curvature $\hat H^t$ does not depend on $t$ since $\hat M\cap\hat\ga$ is a 
parallel dodecagon, and denote  
$$
\hat H^t=\hat H=m_1(\hat \la+\hat \nu)+m_2(\hat \mu+\hat \tau).  
$$
On the other hand, the mean curvature $H^t$ of $M$ is obtained from (\ref{hatla4}) 
$$
H^t=a_t\hat H^t+c_tK,\quad t=1,\dots,8, 
$$
 where $K=2(m_1+m_2)$. 
 Because $H^t$ is constant, we have 
$$
0=H^1-H^2=(a_1-a_2)\hat H+(c_1-c_2)K=2x(\hat u \hat H-\hat v K) 
=2x\hat u(\hat H-\hat\la K).
$$
Since $\hat\la$ is the largest principal curvature, we have $\hat H-\hat\la K<0$, 
and it follows 
\beq
x=0.
\label{x0}
\enq
Next, from 
$$
\bary
0&=H^1-H^3=(a_1-a_3)\hat H+(c_1-c_3)K=2y((\hat v-\hat u)\hat H+(\hat u+\hat v))K \\ 
&=2y(\hat v-\hat u)(\hat H-\hat \tau K)
\eary
$$
where $(\hat v-\hat u)(\hat H-\hat \tau K)>0$ by Lemma \ref{lhatla}  
and as $\hat \tau$ is the smallest principal curvature, we have 
\beq
y=0. 
\label{y0}
\enq
Therefore, we obtain:
\bep
 $C$ is an isometry $\bpm T&0\\0&I'\epm\in   O(2,2)$, $I'=\bpm \eps&0\\0&1\epm$, 
and $M\cap\ga$ itself is a parallel octagon isometric to $M_\theta\cap\ga_\theta$,  
where $\theta$ is uniquely determined from $H$ of $M$ (Fact 4).  
\label{piso}
\enp

\smallskip\noindent
{\em Proof of Theorem \ref{tmain} (iii):} In the same way as the proof of Theorem \ref{tmain} (ii), 
$M\cap\ga$ at a maximum point $p=p^1$ of $\la$ with $\mu$ critical on $L(p^1)$,  
and $M\cap\ga'$ at a minimum point $q=q^1$ of $\la$ with $\mu$ critical on $L(q^1)$, 
are both congruent with a parallel octagon $M_\theta\cap\ga_\theta$,  
where $\theta$ is uniquely determined by the mean curvature of $M$. 
Thus, $\la$ is constant. 
A similar argument at maximum and minimum points of $\mu$ implies 
$\mu$ is constant. Then (\ref{HPhi}) uniquely determines $\nu>\tau$ as 
solutions of a quadratic equation similar to (\ref{AB}), 
and all the principal curvatures are constant all over $M$.
Thus, $M$ is isoparametric. 
\qed

\section{Closed case: $g=6$ (iv)}
\label{sg6}

\subsection{Strategy for $g=6$}
\label{ssg6}

Now we consider the case $g=6$,  
 when the multiplicity  $m=m_1=m_2$ takes values 1 or 2 \cite{A}. 
Let $\la_1> \dots>\la_6 $ be the distinct principal curvatures, 
which  we also denote by $\la=\la_1,\mu=\la_2,\nu=\la_3, \rho=\la_4,\si=\la_5 $ and $\tau=\la_6$, 
using the same notation as in \cite{M6}. 
Later, we use upper indices for points $p^t$.  
The corresponding curvature distributions decompose the tangent space into  $D_1\oplus\cdots\oplus D_6$. 
We denote by $e_j$ any unit vector in $D_j$.

We assume CMC and put 
\beq
\bary
H=m\sum_{i=1}^6\la_i, \quad m=1,2.  
\label{ms6}
\eary
\enq
Note that there are essentially three Lie curvatures when $g=6$ \cite{M6}.   
Now, we prove:

\smallskip\noindent
{\bf Theorem 1.1} (iv)  {\em Let $M$ be a closed CMC Dupin hypersurface with $g=6$.
If $M$ has three independent constant Lie curvatures, 
then $M$ is isoparametric. }
\label{tg6}

\smallskip
\Remark
Even if all Lie curvatures are constant, a Dupin hypersurface with $g=6$ is not 
necessarily Lie equivalent to an isoparametric hypersurface \cite{M6}.

Recall for any normal geodesic $\bar\ga$ of the minimal isoparametric hypersurface $\bar M$, 
$\bar M\cap\bar\ga$ is a regular dodecagon. 
In \cite{M6} we prove: 

\noindent
{\bf Fact 6. (Lemma 4 in \cite{M6}) }
{\it When all the Lie curvatures of a Dupin hypersurface $M$ with $g=6$ are constant,  
there exists a Lie transformation at each point of $p\in M$ which maps 
$M\cap\ga$ into a regular dodecagon, where $\ga$ is the normal geodesic of $M$ at $p$.  }
\label{fLietr}

\smallskip
As before, this is a local fact, and no global result comes out from this,  
but since Lie curvatures are Lie invariant,  
the values are computed from the principal curvatures of 
$\bar M$:  
i.e., putting $g=6$ and $\theta=0$ in  (\ref{thetai}):  
\beq
\la=2+\sqrt3=-\tau,\quad \mu=1=-\si,\quad \nu=2-\sqrt3=-\rho. 
\label{pcmini}
\enq

Now, using CMC and the explicit values of the Lie curvatures, 
we find a point $p\in M$ at which all $d_{ji}$ vanish (Proposition \ref{pdji0}). 
Our goal is to show $M\cap\ga$ at this point is {\em itself} a parallel dodecagon (Proposition \ref{ppara6}). 
Then as in the case $g=4$, the mean curvature uniquely determines $\theta$ so that $M\cap \ga$ is 
isometric to $M_\theta\cap\ga_\theta$, and the maximum and minimum 
values of $\la$ turn out to coincide. Then $\la$ is constant all over $M$,  
and by further argument, we show  $M$ is isoparametric. 

\subsection{Critical point of principal curvatures}
\label{ssCritical}

We use indices $1\le i, j\le 6$. 
Then from (\ref{ms6}),  it follows for $e_j\in D_j $,
\beq
d_{j1}+d_{j2}+d_{j3}+d_{j4}+d_{j5}+d_{j6}=0, \quad j=1,\dots,6. 
\label{dijpf}
\enq
Since the multiplicities of $\la_i$'s are common $m=1,2$, we can omit it.  
As before, it is essential to find a point at which all $d_{ji}=0$ 
under the conditions CMC and constant Lie curvatures (Proposition \ref{pdji0}).  
This is the longest part of the proof. 

For instance, in the case $m=1$, 
there are 30 unknowns $d_{ji}=e_j(\la_i)$, $1\le i, j\le 6$, since $d_{jj}=0$ as $M$ is Dupin.
Because (\ref{dijpf}) consists of 6 equations, the unknowns reduce to 24, in general. 
Take a critical point $p\in M$ of $\la$, and further $\mu$ is critical at $p$ on $L_1(p)$. 
Then we have at $p$, 
\beq
d_{j1}=0, \quad d_{12}=0, \quad 1\le j \le 6. 
\label{dj10}
\enq
These are 6 equations as $d_{11}=0$ already, and the unknown is reduced to 18. 
Now consider three Lie curvatures 
\beq
\Phi_h=\dfrac{(\la-\mu)(\la_h-\si  )}{(\la-\si  )(\la_h-\mu)},\quad h=3,4,6, 
\label{Phih}
\enq
and assume that each $\Phi_h$ is constant on $M$ for $h=3,4,6$. 
Then from
\beq
e_j(\Phi_h)=0,\quad j=1,\dots,6,\quad  h=3,4,6, 
\label{56}
\enq
we have 18 more equations among $d_{ji}$'s for 18 unknowns, and  
we can discuss if $d_{ji}(p)=0$ holds or not.

\medskip
In the following, we investigate each process for $m=1,2$. 
\bel
At $p\in M$ where $\la$ is critical on $M$ and $\mu$ is critical on $L_1(p)$, 
all $d_{1i}$  vanish for any  $e_1\in D_1(p)$ and $i\in\{1,\dots,6\}$. 
\label{ld1i0}
\enl
\Proof
From (\ref{Phih}), we have  for $h\in \{3, 4,6\}$,  
\beq
\bary
e_j(\log\Phi_h)&=\dfrac{d_{j1}-d_{j2}}{\la-\mu}
-\dfrac{d_{j1}-d_{j5}}{\la-\si  }+\dfrac{d_{jh}-d_{j5}}{\la_h-\si  }
-\dfrac{d_{jh}-d_{j2}}{\la_h-\mu}\\
&~\\
&=u_{h} d_{j2}+ v_{h} d_{j5}+ w_{h} d_{jh}=0,\quad j=1,\dots, 6, \,\, h=3, 4, 6, 
\eary
\label{albe}
\enq
where 
\beq
\bary
 u_{h} &=\dfrac1{\la_h-\mu}-\dfrac1{\la-\mu}=\dfrac{(\la-\la_h)}{(\la_h-\mu)(\la-\mu)}\\
 v_{h}&=\dfrac1{\la-\si }-\dfrac1{\la_h-\si}=\dfrac{(\la_h-\la)}{(\la-\si )(\la_h-\si )}\\
 w_{h} & =\dfrac1{\la_h-\si}-\dfrac1{\la_h-\mu}=\dfrac{(\si-\mu)}{(\la_h-\si )(\la_h-\mu)}. 
 \eary
\label{434}
\enq
By the assumption,  $d_{j1}=0$ and $d_{12}=0$ hold at $p$. 
Put $j=1$ in (\ref{albe}) and using $d_{12}=0$ at $p$, we have 
\beq
\bary 
v_3d_{15}+w_3d_{13}=0,\quad \text{i.e.,} \quad d_{13}=-\dfrac{v_3}{w_3}d_{15},\\ 
v_4d_{15}+w_4d_{14}=0,\quad \text{i.e.,}\quad d_{14}=-\dfrac{v_4}{w_4}d_{15},\\ 
v_6d_{15}+w_6d_{16}=0,\quad \text{i.e.,}\quad d_{16}=-\dfrac{v_6}{w_6}d_{15}. 
\eary
\label{d1h}
\enq
Substituting these into (\ref{dijpf}), we obtain 
\beq
\Big(1-\dfrac{v_3}{w_3}-\dfrac{v_4}{w_4}-\dfrac{v_6}{w_6}\Big)d_{15}=0.
\label{d1i}
\enq
Since $v_3,v_4<0$,  $v_6>0$,   $w_3,w_4>0$ and $w_6<0$, we have 
$$
1-\dfrac{v_3}{w_3}-\dfrac{v_4}{w_4}-\dfrac{v_6}{w_6}>0,  
$$
and it follows 
 \beq
 d_{1i}=0.
 \label{d1i0}
 \enq
 \qed
 
 \bel
On $L_1(p)$ where $\la$ is critical, 
 $d_{2i}=d_{3i}=d_{4i}=d_{6i}=0$. 
\label{ldji06}
\enl
\Proof
We do not use $d_{12}=0$ in the proof below. Hence the argument holds all over $L_1(p)$. 
\begin{enumerate}
\item First, to show $d_{2i}=0$ on  $L_1(p)$, 
put $j=2$ in (\ref{albe}), and we have 
\beq
\bary
 v_{3} d_{25}+ w_{3} d_{23}=0,\quad \text{i.e.,}\quad d_{23}=-\dfrac{v_3}{w_3}d_{25}\\
  v_{4} d_{25}+ w_{4} d_{24}=0,\quad  \text{i.e.,}\quad d_{24}=-\dfrac{v_4}{w_4}d_{25}\\
   v_{6} d_{25}+ w_{6} d_{26}=0,\quad  \text{i.e.,}\quad d_{26}=-\dfrac{v_6}{w_6}d_{25}.  
   \eary
  \label{d23}
  \enq
These are same as (\ref{d1h}) replacing $d_{1j}$ with $d_{2j}$. 
Thus we obtain 
  \beq
 d_{2i}=0.
 \label{d2i0}
 \enq

\item To obtain $d_{3i}=0$, consider 
  \beq
 \check\Psi_h=\dfrac{(\nu-\rho)(\la_h-\la)}{(\nu-\la)(\la_h-\rho)},\quad h=2,5,6. 
\enq
Then we have, using $d_{31}=0$ and $d_{33}=0$, 
\beq
\bary
0&=e_3(\log \check\Psi_h)=\dfrac{-d_{34}}{\nu-\rho}+\dfrac{d_{3h}}{\la_h-\la}-\dfrac{d_{3h}-d_{34}}{\la_h-\rho}\\
&=\Big(\dfrac1{\la_h-\la}-\dfrac1{\la_h-\rho}\Big)d_{3h}+\Big(\dfrac1{\la_h-\rho}-\dfrac1{\nu-\rho}\Big)d_{34}\\
&=\dfrac1{\la_h-\rho}\Big(\dfrac{\la-\rho}{\la_h-\la}d_{3h}+\dfrac{\nu-\la_h}{\nu-\rho}d_{34}\Big).
\eary
\label{d3h}
\enq
Putting $h=2,5,6$, we obtain
\beq
\bary
d_{32}=\dfrac{(\la-\mu)(\nu-\mu)}{(\la-\rho)(\nu-\rho)}d_{34},\quad
d_{35}=\dfrac{(\la-\si)(\nu-\si)}{(\la-\rho)(\nu-\rho)}d_{34},\\
d_{36}=\dfrac{(\la-\tau)(\nu-\tau)}{(\la-\rho)(\nu-\rho)}d_{34}. 
\eary
\label{}
\enq
Here, we use the Lie curvature.  Using (\ref{pcmini}), we compute 
\beq
\bary
\dfrac{(\la-\mu)(\nu-\mu)}{(\la-\rho)(\nu-\rho)}=\dfrac{(\la-\rho)(\nu-\mu)}{(\la-\mu)(\nu-\rho)}\Big(\dfrac{\la-\mu}{\la-\rho}\Big)^2
=-2\Big(\dfrac{\la-\mu}{\la-\rho}\Big)^2,\\
\dfrac{(\la-\si)(\nu-\si)}{(\la-\rho)(\nu-\rho)}
=\dfrac{(\la-\rho)(\nu-\si)}{(\la-\si)(\nu-\rho)}\Big(\dfrac{\la-\si}{\la-\rho}\Big)^2
=2\Big(\dfrac{\la-\si}{\la-\rho}\Big)^2\\
\dfrac{(\la-\tau)(\nu-\tau)}{(\la-\rho)(\nu-\rho)}
=\dfrac{(\la-\rho)(\nu-\tau)}{(\la-\tau)(\nu-\rho)}\Big(\dfrac{\la-\tau}{\la-\rho}\Big)^2
= 4\Big(\dfrac{\la-\tau}{\la-\rho}\Big)^2.
\eary
\label{}
\enq
Therefore, we have 
\beq
\bary
0&=d_{32}+d_{34}+d_{35}+d_{36}\\
&=d_{34}\Big(-2\Big(\dfrac{\la-\mu}{\la-\rho}\Big)^2+1+2\Big(\dfrac{\la-\si}{\la-\rho}\Big)^2
+4\Big(\dfrac{\la-\tau}{\la-\rho}\Big)^2\Big)\\
\eary
\enq
but as the coefficient satisfies
$$
-2\Big(\dfrac{\la-\mu}{\la-\rho}\Big)^2+1+2\Big(\dfrac{\la-\si}{\la-\rho}\Big)^2
+4\Big(\dfrac{\la-\tau}{\la-\rho}\Big)^2>-2+1+2+4=5>0,
$$
we obtain 
\beq
d_{3i}=0.
\label{d3i0}
\enq

 \item To obtain $d_{4i}=0$, consider 
\beq
\bar\Psi_h=\dfrac{(\rho-\nu)(\la_h-\la)}{(\rho-\la)(\la_h-\nu)},\quad h=2,5,6. 
\enq
Then we have, using $d_{41}=0$ and $d_{44}=0$, 
\beq
\bary
0&=e_4(\log\bar\Psi_h)=\dfrac{-d_{43}}{\rho-\nu}+\dfrac{d_{4h}}{\la_h-\la}-\dfrac{d_{4h}-d_{43}}{\la_h-\nu}\\
&=\Big(\dfrac1{\la_h-\la}-\dfrac1{\la_h-\nu}\Big)d_{4h}+\Big(\dfrac1{\la_h-\nu}-\dfrac1{\rho-\nu}\Big)d_{43}\\
&=\dfrac1{\la_h-\nu}\Big(\dfrac{\la-\nu}{\la_h-\la}d_{4h}+\dfrac{\rho-\la_h}{\rho-\nu}d_{43}\Big).
\eary
\label{d3h}
\enq
Putting $h=2,5,6$, we obtain
\beq
\bary
d_{42}=\dfrac{(\la-\mu)(\rho-\mu)}{(\la-\nu)(\rho-\nu)}d_{43},\quad
d_{45}=\dfrac{(\la-\si)(\rho-\si)}{(\la-\nu)(\rho-\nu)}d_{43},\\
d_{46}=\dfrac{(\la-\tau)(\rho-\tau)}{(\la-\nu)(\rho-\nu)}d_{43}
\eary
\label{}
\enq
Using the Lie curvatures which we compute from (\ref{pcmini}), we have
\beq
\bary
\dfrac{(\la-\mu)(\rho-\mu)}{(\la-\nu)(\rho-\nu)}
=\dfrac{(\la-\nu)(\rho-\mu)}{(\la-\mu)(\rho-\nu)}\Big(\dfrac{\la-\mu}{\la-\nu}\Big)^2
=3\Big(\dfrac{\la-\mu}{\la-\rho}\Big)^2,\\
\dfrac{(\la-\si)(\rho-\si)}{(\la-\nu)(\rho-\nu)}
=\dfrac{(\la-\nu)(\rho-\si)}{(\la-\si)(\rho-\nu)}\Big(\dfrac{\la-\si}{\la-\nu}\Big)^2
=-\Big(\dfrac{\la-\si}{\la-\rho}\Big)^2\\
\dfrac{(\la-\tau)(\rho-\tau)}{(\la-\nu)(\rho-\nu)}
=\dfrac{(\la-\nu)(\rho-\tau)}{(\la-\tau)(\rho-\nu)}\Big(\dfrac{\la-\tau}{\la-\nu}\Big)^2
= -3\Big(\dfrac{\la-\tau}{\la-\nu}\Big)^2. 
\eary
\label{}
\enq
Then we obtain 
\beq
\bary
0=d_{42}+d_{43}+d_{45}+d_{46}\\=
d_{43}\Big(3\Big(\dfrac{\la-\mu}{\la-\rho}\Big)^2+1-\Big(\dfrac{\la-\si}{\la-\rho}\Big)^2
 -3\Big(\dfrac{\la-\tau}{\la-\nu}\Big)^2\Big)\\
\eary
\enq
where the coefficient satisfies
$$
3\Big(\dfrac{\la-\mu}{\la-\rho}\Big)^2+1-\Big(\dfrac{\la-\si}{\la-\rho}\Big)^2
 -3\Big(\dfrac{\la-\tau}{\la-\nu}\Big)^2<3+1-1-3=0.
$$
Thus we have  
\beq
d_{4i}=0.
\label{d4i0}
\enq
 
  \item For $d_{6i}$, consider 
 \beq
\tilde\Psi_h=\dfrac{(\tau-\nu)(\la_h-\la)}{(\tau-\la)(\la_h-\nu)},\quad h=2,4,5. 
\enq  
Then we have, using $d_{61}=0$ and $d_{66}=0$, 
$$
\bary
0&=e_6(\log \tilde\Psi_h)=
\dfrac{-d_{63}}{\tau-\nu}+\dfrac{d_{6h}}{\la_h-\la}-\dfrac{d_{6h}-d_{63}}{\la_h-\nu}\\
&=\Big(\dfrac1{\la_h-\la}-\dfrac1{\la_h-\nu}\Big)d_{6h}+\Big(\dfrac1{\la_h-\nu}-\dfrac1{\tau-\nu}\Big)d_{63}\\
&=\dfrac1{\la_h-\nu}\Big(\dfrac{\la-\nu}{\la_h-\la}d_{6h}+\dfrac{\tau-\la_h}{\tau-\nu}d_{63}\Big),
\eary
$$
which implies
\beq
\bary
d_{62}=\dfrac{(\la-\mu)(\tau-\mu)}{(\la-\nu)(\tau-\nu)}d_{63},\quad
d_{64}=\dfrac{(\la-\rho)(\tau-\rho)}{(\la-\nu)(\tau-\nu)}d_{63},\\
d_{65}=\dfrac{(\la-\si)(\tau-\si)}{(\la-\nu)(\tau-\nu)}d_{63}.
\eary
\label{d6i}
\enq
Then the coefficient of $d_{63}$ in $d_{61}+\dots+d_{65}$ is positive, and   we have  
   \beq
 d_{6i}=0.
 \label{d6i0}
 \enq
 \end{enumerate}
 \qed

 To see if $d_{5i}=0$, consider
  \beq
\tilde\Psi_h=\dfrac{(\si-\mu)(\la_h-\la)}{(\si-\la)(\la_h-\mu)},\quad h=3,4,6. 
\enq  
Then using $d_{51}=0$ and $d_{55}=0$, we have 
$$
\bary
0&=e_5(\log \tilde\Psi_h)=
\dfrac{-d_{52}}{\si-\mu}+\dfrac{d_{5h}}{\la_h-\la}-\dfrac{d_{5h}-d_{52}}{\la_h-\mu}\\
&=\Big(\dfrac1{\la_h-\la}-\dfrac1{\la_h-\mu}\Big)d_{5h}+\Big(\dfrac1{\la_h-\mu}-\dfrac1{\si-\mu}\Big)d_{52}\\
&=\dfrac1{\la_h-\mu}\Big(\dfrac{\la-\mu}{\la_h-\la}d_{5h}+\dfrac{\si-\la_h}{\si-\mu}d_{52}\Big),
\eary 
$$
which implies
\beq
\bary
d_{53}=\dfrac{(\la-\nu)(\si-\nu)}{(\la-\mu)(\si-\mu)}d_{52},\quad
d_{54}=\dfrac{(\la-\rho)(\si-\rho)}{(\la-\mu)(\si-\mu)}d_{52},\\
d_{56}=\dfrac{(\la-\tau)(\si-\tau)}{(\la-\mu)(\si-\mu)}d_{52}.
\eary
\label{d5h}
\enq
Thus we have
$$
0=d_{52}\Big(1+\dfrac{(\la-\nu)(\si-\nu)}{(\la-\mu)(\si-\mu)}+\dfrac{(\la-\rho)(\si-\rho)}{(\la-\mu)(\si-\mu)}
+\dfrac{(\la-\tau)(\si-\tau)}{(\la-\mu)(\si-\mu)}\Big),
 $$
but we cannot determine the signature of the coefficient as before even if we use the values of Lie curvatures. 
Suppose $d_{52}\ne0$ occurs on an open neighborhood $U$ on $L_1$. 
 Then we have 
 \beq
(\la-\mu)(\si-\mu)+(\la-\nu)(\si-\nu)+ (\la-\rho)(\si-\rho)+(\la-\tau)(\si-\tau)=0, 
\label{d5i1}
\enq
i.e., 
\beq
\mu^2+\nu^2+\rho^2+\tau^2-(\la+\si)(\mu+\nu+\rho+\tau)+4\la\si=0 
\label{d5i}
\enq
on $U$. Since $d_{51}=d_{55}=0$, and $e_5(H)=0$, we have
$$
\bary
0=e_5\big(\mu^2+\nu^2+\rho^2+\tau^2-(\la+\si)(\mu+\nu+\rho+\tau)+4\la\si\big)\\
=2(\mu d_{52}+\nu d_{53}+\rho d_{54}+\tau d_{56})-(\la+\si)(d_{52}+d_{53}+d_{54}+d_{56})\\
=2(\mu d_{52}+\nu d_{53}+\rho d_{54}+\tau d_{56}), 
\eary
$$
and using (\ref{d5h}) again, we have 
$$
0=d_{52}\Big(\mu+\nu\dfrac{(\la-\nu)(\si-\nu)}{(\la-\mu)(\si-\mu)}+\rho\dfrac{(\la-\rho)(\si-\rho)}{(\la-\mu)(\si-\mu)}
  +\tau\dfrac{(\la-\tau)(\si-\tau)}{(\la-\mu)(\si-\mu)}\Big). 
$$  
This implies, as $d_{52}\ne0$ on $U$, 
\beq
\mu(\la-\mu)(\si-\mu)+\nu(\la-\nu)(\si-\nu)+ \rho(\la-\rho)(\si-\rho)+\tau(\la-\tau)(\si-\tau)=0. 
\label{mulamu}
\enq
Because (\ref{d5i1})  implies
$$
\tau(\la-\tau)(\si-\tau)=-\tau\big((\la-\mu)(\si-\mu)+(\la-\nu)(\si-\nu)+ (\la-\rho)(\si-\rho)\big),
$$
substituting this into (\ref{mulamu}), we obtain 
$$
(\mu-\tau)(\la-\mu)(\si-\mu)+(\nu-\tau)(\la-\nu)(\si-\nu)+ (\rho-\tau)(\la-\rho)(\si-\rho)=0.
$$
However, each term of the LHS is negative, a contradiction.
Therefore, we know  
\beq
d_{5i}=0. 
\label{d5i0}
\enq
As a conclusion:
\bep All $d_{ji}$ vanish at $p$.  
\label{pdji0}
\enp

\subsection{Angle relation and the Lie curvature}
\label{ssAngle}

\begin{figure}[htbp]
\centering
\begin{minipage}[b]{0.50\columnwidth}
    \centering
    \includegraphics[width=0.9\columnwidth]{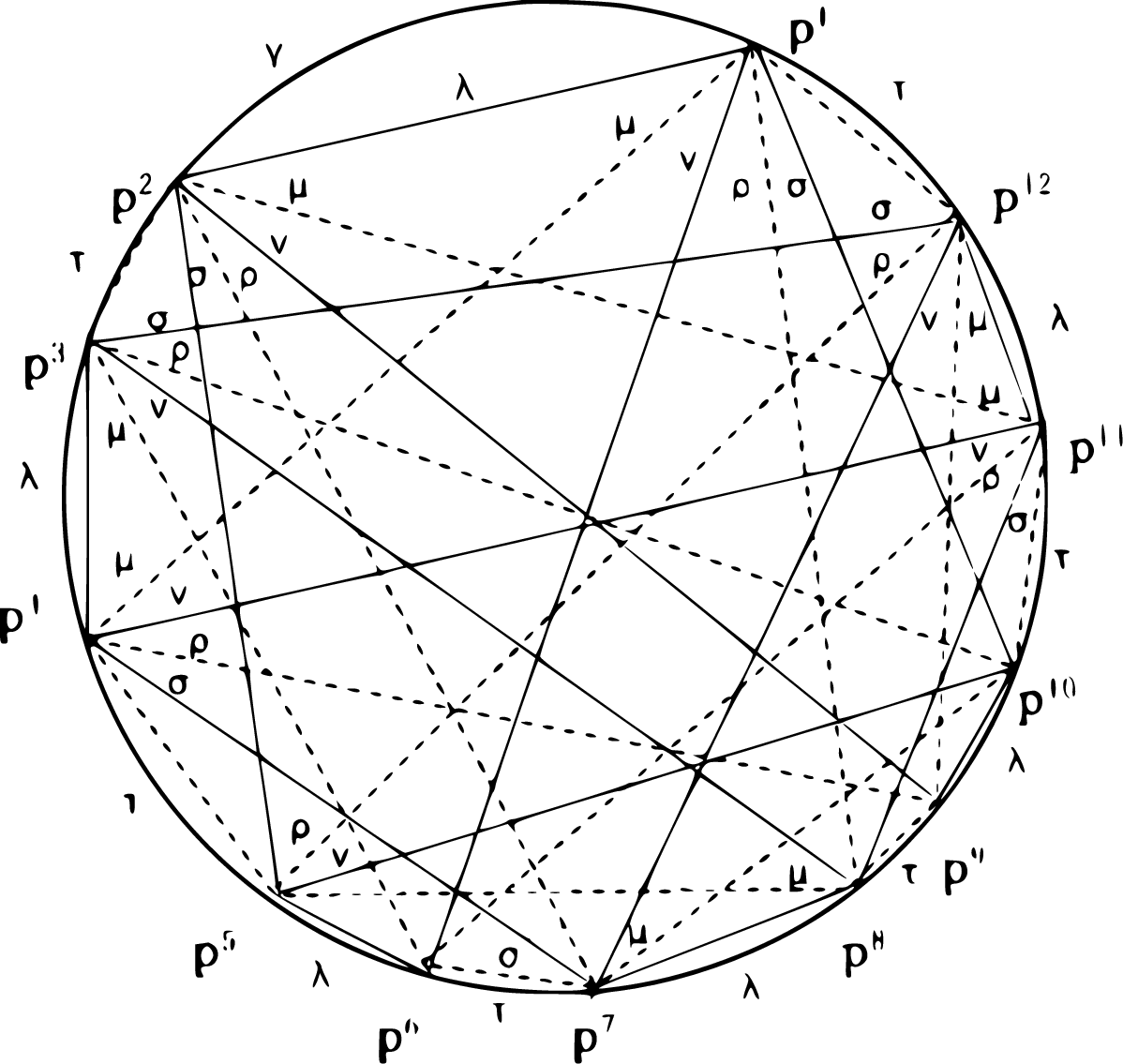}
    \caption{Link of leaves}
    \label{flinkg6}
\end{minipage}
\begin{minipage}[b]{0.49\columnwidth}
    \centering
    \includegraphics[width=59mm]{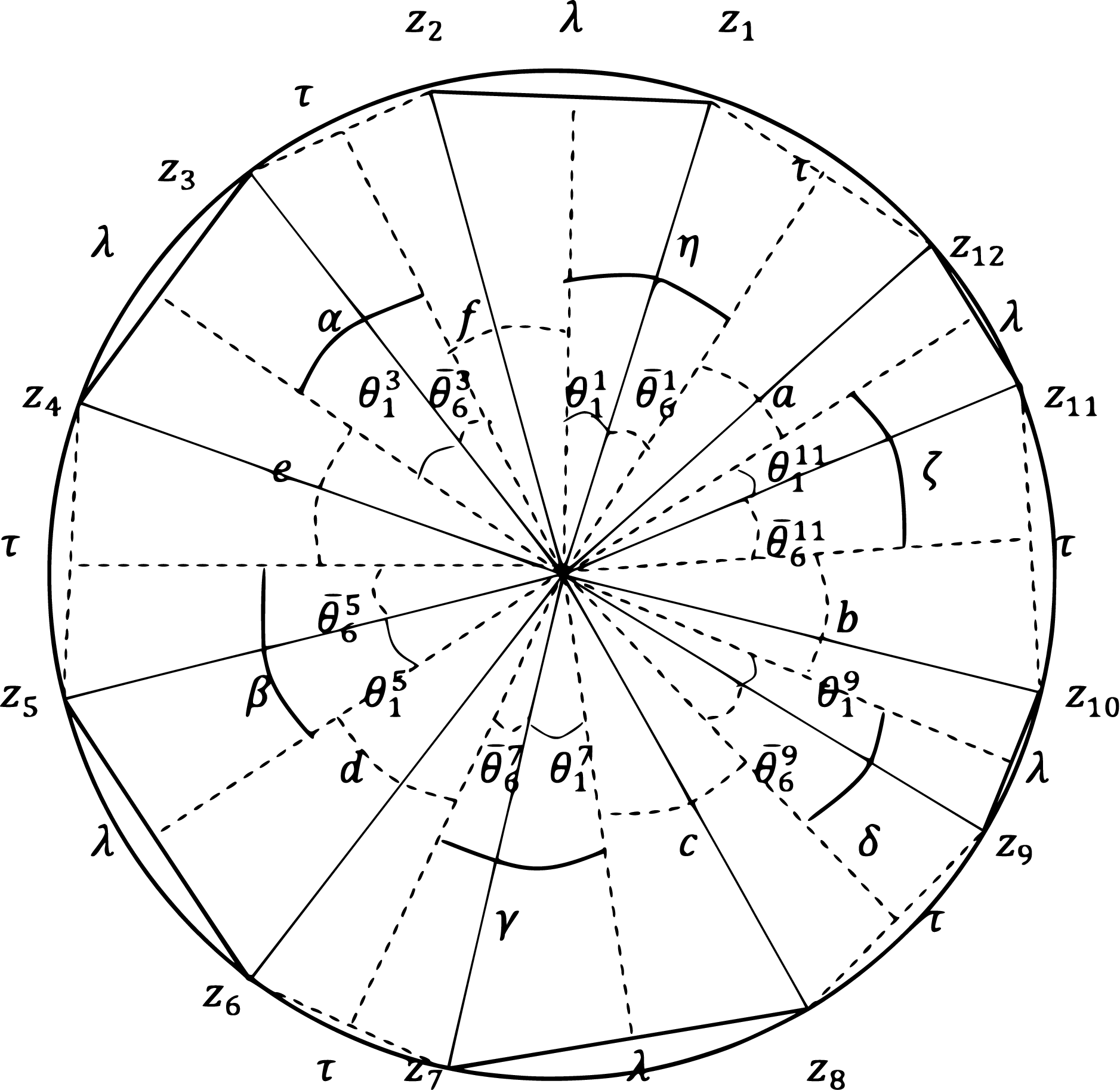}
\caption{Angle relation}
\label{fdodec}
\end{minipage}
\end{figure}

In \cite{M6}, we show that if  all $d_{ji}$ vanish at  $p$,  
all the leaves $L_i$ through $p^1=p$ are  
totally geodesic in the curvature sphere $C_i$,  and $L_i\cap \gamma=\{p^1,p^{2i}\}$, 
where $\ga$ is the normal geodesic through $p^1$, and $p^{2i}$ is  the antipodal point of $p^1$ in $L_i$.  
Now from  Lemma 3 in \cite{M6} where we use the tautness, it follows  
\bel
We have $M \cap \gamma=\{p=p^1,\dots,p^{12}\}$, 
where $p^1, p^{2i}$, $i=2,3,4,5,6$ are as above, and $p^{\text{odd}}$ are as in Figure \ref{flinkg6}. 
The leaves at each $p^t$ are all totally geodesic 
in the curvature sphere and cut $\gamma$ orthogonally at some $p^s$. 
They intersect with each other, as is shown in Figure \ref{flinkg6}. 
\label{lint}
\enl

For an isoparametric hypersurface $M_\theta$,  
$M_\theta\cap \ga_\theta$ is a parallel dodecagon (Figure \ref{fparag6})  for any normal geodesic 
$\ga_\theta$ of $M_\theta$.   
The goal is to show that $M\cap\ga=\{p^1,\dots,p^{12}\}$ is itself a parallel dodecagon.  
We prove it in several steps.  

\bep
Denote $\la_i^t=\la_i(p^t)=\cot \theta_i^t$, $0<\theta_i^t<\pi$,  for $i=1,\dots,6$, and $t=1,\dots,12$.  
Putting 
\beq
\bary
\theta^1_2=\theta^1_1+\al,\, \theta^1_3=\theta^1_2+ \be,\, \theta^1_4=\theta^1_3 +\ga,
\theta^1_5=\theta^1_4+\de,\, \theta^1_6=\theta^1_5+\zeta\\
\theta^2_2=\theta^2_1+a,\, \theta^2_3=\theta^2_2 +b,\, \theta^2_4=\theta^2_3+c,
\theta^2_5=\theta^2_4+d,\, \theta^2_6=\theta^2_5+e, 
\eary
\enq
and  $\eta=\pi-(\al+\be+\ga+\de+\zeta)$, $f=\pi-a-b-c-d-e$, 
we can express the angles between $\theta^t_i$ and $\theta^t_{i+1}$ at $p^t$, $t=1,\dots,12$   
as follows: 

\medskip
{\small
\begin{tabular}{|l|r|r|r|r|r|r|r|r|r|r|r|r|r|} \hline
$p^t$ & $\theta_1^t$ &  & $\theta_2^t$ & & $\theta_3^t$ & & $\theta_4^t $ & & $\theta^t_5$ && $\theta^t_6$ && 
$\pi+\theta^t_1$\\ \hline 
 $p^1$ &  &$\al$ & & $\be$ && $\ga$ &&$\de$ && $\zeta$ && $\eta$ &\\
 $p^2$ &  &$a$ & & $b$ && $c$ &&$d$ && $e$ && $f$ &\\
$p^3$ & &$\be$ && $\ga$ && $\de$ && $\zeta$ && $\eta$ && $\al$ &\\
$p^4$ & &$f$ && $a$ && $b$ && $c$ && $d$ && $e$ &\\
$p^5$ & &$\ga$ && $\de$ && $\zeta$ && $\eta$ && $\al$ && $\be$ &\\
$p^6$ & &$e$ && $f$ && $a$ && $b$ && $c$ && $d$ &\\
$p^7$ & &$\de$ && $\zeta$ && $\eta$ && $\al$ && $\be$ && $\ga$ &\\
$p^8$ & &$d$ && $e$ && $f$ && $a$ && $b$ && $c$ &\\
$p^9$ & & $\zeta$ && $\eta$ && $\al$ && $\be$ && $\ga$ &&$\de$ &\\
$p^{10}$ & & $c$ && $d$ && $e$ && $f$ && $a$ &&$b$ &\\
$p^{11}$ & &$\eta$ && $\al$ && $\be$ && $\ga$ && $\de$ && $\zeta$ &\\
$p^{12}$ & &$b$ && $c$ && $d$ && $e$ && $f$ && $a$ &\\
 \hline
 \end{tabular}
  \captionof{table}{Angles 2} 
}
\label{palbega6}
\enp

\Proof
We apply a similar argument as Proposition \ref{palbega4} and describe it shortly. 
Denote $p^t$ by the complex number $z^t$, and put  
$\bar\theta_i^t=\pi-\theta_i^t$, $0< \bar\theta_i^t<\pi$.  
As before, from Figure \ref{fdodec}, taking the orientation into account, we have
$$
\theta^1_2=\theta^1_1+\theta^3_1+\bar\theta^3_6=\theta^1_1+\al,\quad 
\theta^1_3=\theta^1_2+\theta^5_1+\bar\theta^5_6=\theta^1_1+\be,\dots
$$
In this way, we obtain the table for $p^{\text{odd}}$, and starting from $p^2=z_2$ clock-wise, 
we obtain the table for $p^{\text{even}}$.
\qed

\smallskip
In the following, all the angles stay in $(0,\pi)$ modulo $\pi$, and we can use  
(\ref{LPhi}), where taking care of the order of the principal curvatures, we have  
\beq
\Psi_\nu(z_1)=\dfrac{(\la-\mu)(\nu-\si)}{(\la-\si)(\nu-\mu)}=[z_2,z_6;z_4,z_{10}]=-1.
\label{Psiz}
\enq
Here, the value of $\Psi$ is computed from (\ref{pcmini}).

\smallskip
For the moment,  we advance our argument under the assumption $\la^1\rho^1=-1$ at $p^1\in M$. 
In \S \ref{ssconformal}, we apply this argument to the conformal image of $M$ satisfying this assumption. 
\bel
Assume $\la^1\rho^1=-1$ at $p^1\in M$.  
Then putting $\Psi_\nu^t=\Psi_\nu(p^t)$ and 
\beq
w_1=e^{2i\al},\, w_2=e^{2i\be},\, w_3=e^{2i\ga},\, w_4=e^{2i\de},\, 
w_5=e^{2i\zeta},\, w_6=e^{2i\eta},
\label{e2i}
\enq
we have
\beq
\bary 
\Psi_\nu^1=[z_2,z_6;z_4,z_{10}]=\dfrac{(1-w_1)(1-w_3w_4)}{(1+w_4)(1+w_1w_3)},\\
\Psi_\nu^5=[z_6,z_{10};z_8,z_{2}]=\dfrac{(1-w_3)(1-w_5w_6)}{(1+w_3)(1+w_5w_6)},\\ 
\Psi_\nu^{11}=[z_{12},z_{4};z_{2},z_{8}]=\dfrac{(1-w_6)(1-w_2w_3)}{(1+w_6)(1+w_2w_3)}.
\eary
\label{Psinu}
\enq
\label{lPsi}
\enl
\Proof
By our assumption,  $2(\al+\be+\ga)=\pi=2(\de+\zeta+\eta)$ holds, and $e^{i(\pi+\vp)}=-e^{i\vp}$. 
Thus from 
 $$
 \bary
 z_4=e^{2i\al}z_2,\quad 
 z_{10}=e^{2i(\al+\be+\ga+\de)}z_2=-z^{2i\de}z_2,\\
  z_{10}=e^{2i(\ga+\de)} z_6,\quad z_4=e^{2i(\ga+\de+\zeta+\eta+\al)}z_6=-e^{2i(\al+\ga)}z_6
 \eary
 $$
  we obtain
 \beq
 \bary 
 \Psi_\nu^1&=[z_2,z_6;z_4,z_{10}]=\dfrac{(z_2-z_4)(z_6-z_{10})}{(z_2-z_{10})(z_6-z_4)}\\
 &=\dfrac{(1-e^{2i\al})(1-e^{2i(\ga+\de)})}{(1+e^{2i\de})(1+e^{2i(\al+\ga)})}
 =\dfrac{(1-w_1)(1-w_3w_4)}{(1+w_4)(1+w_1w_3)}. 
 \eary
 \label{Psi1}
 \enq
 Similarly, we have 
$$
\bary
z_8=e^{2i\ga}z_6,\quad z_2=e^{2i(\ga+\de+\zeta+\eta)}z_6=-e^{2i\ga}z_6\\
z_2=e^{2i(\zeta+\eta)}z_{10},\quad  z_8=e^{2i(\zeta+\eta+\al+\be+\ga)}z_{10}=-e^{2i(\zeta+\eta)}z_{10},
\eary
$$
and we obtain  
$$
\bary
\Psi_\nu^5&=[z_6,z_{10};z_8,z_{2}]=\dfrac{(z_6-z_8)(z_{10}-z_{2})}{(z_6-z_{2})(z_{10}-z_8)}\\
&=\dfrac{(1-e^{2i\ga})(1-e^{2i(\zeta+\eta)})}{(1+e^{2i\ga})(1+e^{2i(\zeta+\eta)})}
=\dfrac{(1-w_3)(1-w_5w_6)}{(1+w_3)(1+w_5w_6)}. 
\eary
$$
Also from 
$$
\bary
z_2=e^{2i\eta}z_{12},\quad z_8=e^{2i(\eta+\al+\be+\ga)}z_{12}=-e^{2i\eta}z_{12}\\
z_8=e^{2i(\be+\ga)}z_4,\quad z_2=e^{2i(\be+\ga+\de+\zeta+\eta)}z_4=-e^{2i(\be+\ga)}z_4, 
\eary
$$
we have 
$$
\bary
\Psi_\nu^{11}&=[z_{12},z_{4};z_{2},z_{8}]=\dfrac{(z_{12}-z_{2})(z_{4}-z_8)}{(z_{12}-z_{8})(z_{4}-z_{2})}\\
&=\dfrac{(1-e^{2i\eta})(1-e^{2i(\beta+\ga)})}{(1+e^{2i\eta})(1+e^{2i(\be+\ga)})}
=\dfrac{(1-w_6)(1-w_2w_3)}{(1+w_6)(1+w_2w_3)}.
\eary
$$
\qed

\bep
Under the assumption $\la^1\rho^1=-1$, we have
\beq
\al=\be=\ga=\de=\zeta=\eta=\dfrac{\pi}6. 
\label{aleta}
\enq
\label{paleta}
\enp
\Proof
In (\ref{Psinu}), using $w_4w_5w_6=-1$, we have from $\Psi_\nu^5=-1$, 
$$
\dfrac{(1-w_3)(1-w_5w_6)}{(1+w_3)(1+w_5w_6)}=\dfrac{(1-w_3)(w_4+1)}{(1+w_3)(w_4-1)}=-1. 
$$
This implies
$$
(1-w_3)(w_4+1)+(1+w_3)(w_4-1)=0,  
$$
and we obtain
\beq
w_3=w_4 \quad \text{i.e.,  }\quad \ga=\de. 
\label{gade}
\enq
Similarly from $\Psi_\nu^{11}=-1$, we have
$$
\dfrac{(1-w_6)(1-w_2w_3)}{(1+w_6)(1+w_2w_3)}=\dfrac{(1-w_6)(w_1+1)}{(1+w_6)(w_1-1)}=-1
$$
and obtain
\beq
w_1=w_6 \quad \text{i.e.,  }\quad \al=\eta.
\label{ale}
\enq
Then from $2(\al+\be+\ga)=\pi=2(\de+\zeta+\eta)$,  it follows 
\beq
\be=\zeta \quad \text{i.e.,  }\quad w_2=w_5. 
\label{bez}
\enq
Next from (\ref{Psi1}), 
we have  
$$
\dfrac{(1-w_1)(1-w_3w_4)}{(1+w_4)(1+w_1w_3)}
=\dfrac{(1-w_1)(1-w_3^2)}{(1+w_3)(1+w_1w_3)}
=\dfrac{(1-w_1)(1-w_3)}{(1+w_1w_3)}=-1
$$
and 
\beq
2(w_1w_3+1)=w_1+w_3
\label{w1w3}
\enq
follows.
This equation has two pieces of information on real and imaginary parts, 
and we can solve two unkowns $\al$ and $\ga$. 

Since $0<2\al,2\ga<\pi$, we have $\sin 2\al,\sin 2\ga>0$, and can 
put  $w_1=e^{2i\al}=x+i\sqrt{1-x^2}$,  
and $w_3=e^{2i\ga}=y+i\sqrt{1-y^2}$,  $x,y\in \R$.  
Then (\ref{w1w3}) becomes 
$$
2(x+i\sqrt{1-x^2})(y+i\sqrt{1-y^2})+1)=x+i\sqrt{1-x^2}+y+i\sqrt{1-y^2}, 
$$
which implies  
$$
\begin{cases}
2(xy-\sqrt{1-x^2}\sqrt{1-y^2})+1)=x+y\\
2(x\sqrt{1-y^2})+y\sqrt{1-x^2})=\sqrt{1-x^2}+\sqrt{1-y^2},
\end{cases}
$$
namely
$$
\begin{cases}
(2xy+2-x-y)^2=4 (1-x^2)(1-y^2)\\
(2x-1)^2(1-y^2)=(2y-1)^2(1-x^2). 
\end{cases}
$$
Thus we have 
\beq
\begin{cases}
(x+y)(5(x+y)-4(xy+1))=0\\
(x-y)(5(x+y)-4(xy+1))=0. 
\end{cases}
\label{xy}
\enq
Here, $x\pm y=0$ implies $x=y=0$, and $\al=\ga=\frac{\pi}4$, 
contradicts $0<\be=\frac{\pi}2-(\al+\ga)$. 
If $x+y=0$ and $x-y\ne0$,  we have
$$
0=5(x+y)-4(xy+1)=4(x^2-1)
$$
which implies $x=\cot 2\al=\pm1$ and $\al=0$ or $\frac{\pi}2$, a contradiction. 
If $x+y\ne0$ and $x-y=0$,  we have
$$
0=5(x+y)-4(xy+1)=-4x^2+10x-4=-2(2x-1)(x-2), 
$$
which implies 
$$
x=\cos 2\al=\dfrac12=y=\cos 2\ga
$$
namely,
\beq
\al=\ga=\dfrac{\pi}6,
\label{algapi6}
\enq
and so $\be=\dfrac{\pi}6$ follows, which proves the proposition. 
\qed

When $\la^2\rho^2=-1$, replacing $z^t$ and $w_t$ suitably, 
we apply a similar argument to $a,b,c,d,e,f$ and obtain:
\bep
If $\la^1\rho^1=-1=\la^2\rho^2$ holds,  we have 
\beq
\al=\be=\ga=\de=\zeta=\eta=a=b=c=d=e=f=\dfrac{\pi}6,
\label{pi6}
\enq
and $M\cap\ga$ is a parallel dodecagon. 
\label{ppi6}
\enp

\subsection{Conformal transformation}
\label{ssconformal}

\begin{figure}[htbp]
\centering
\begin{minipage}[b]{0.46\columnwidth}
    \centering
    \includegraphics[width=57mm]{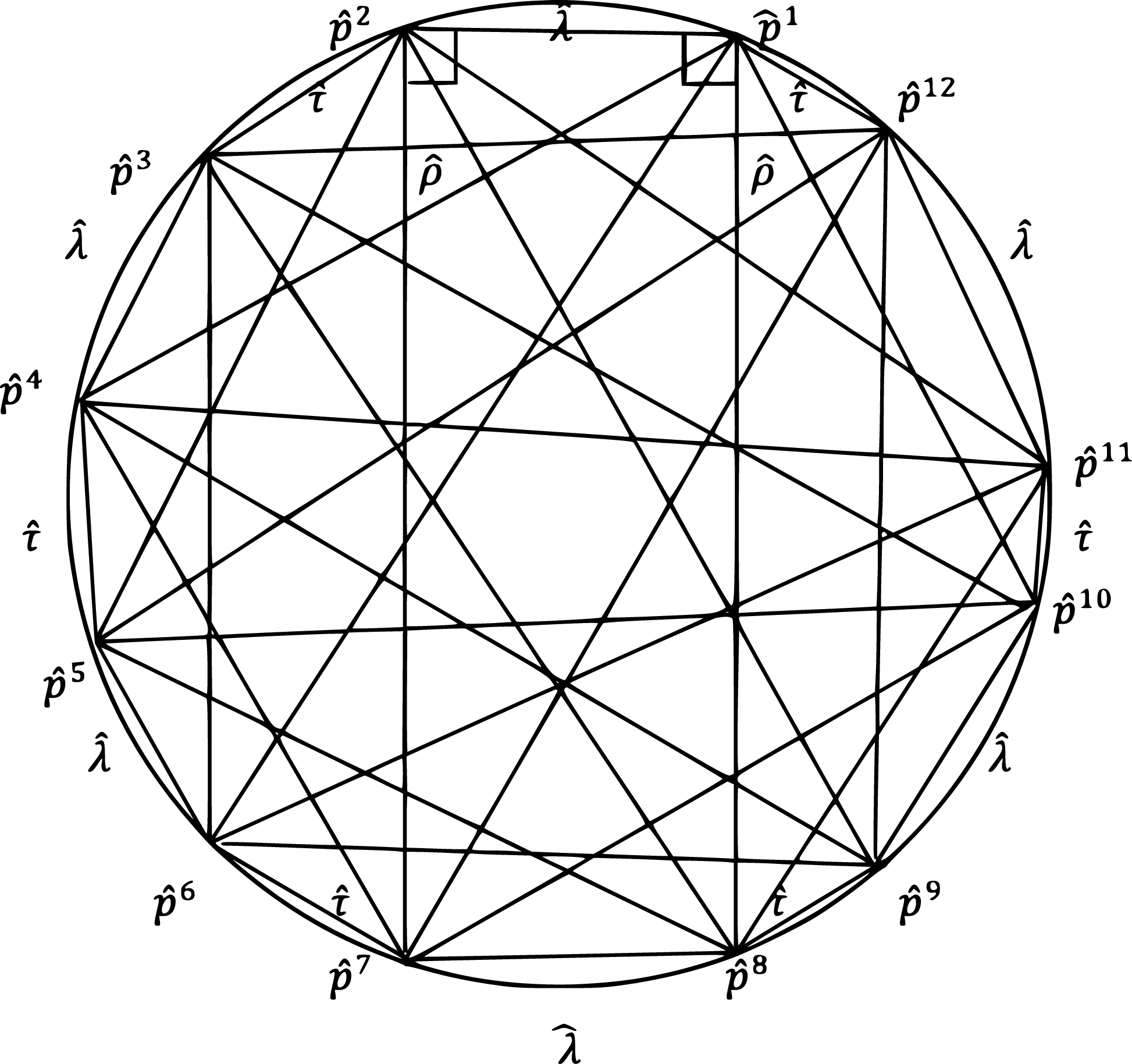}
\caption{Conformal deformation}
\label{fconf}
\end{minipage}
\begin{minipage}[b]{0.45\columnwidth}
    \centering
    \includegraphics[width=55mm]{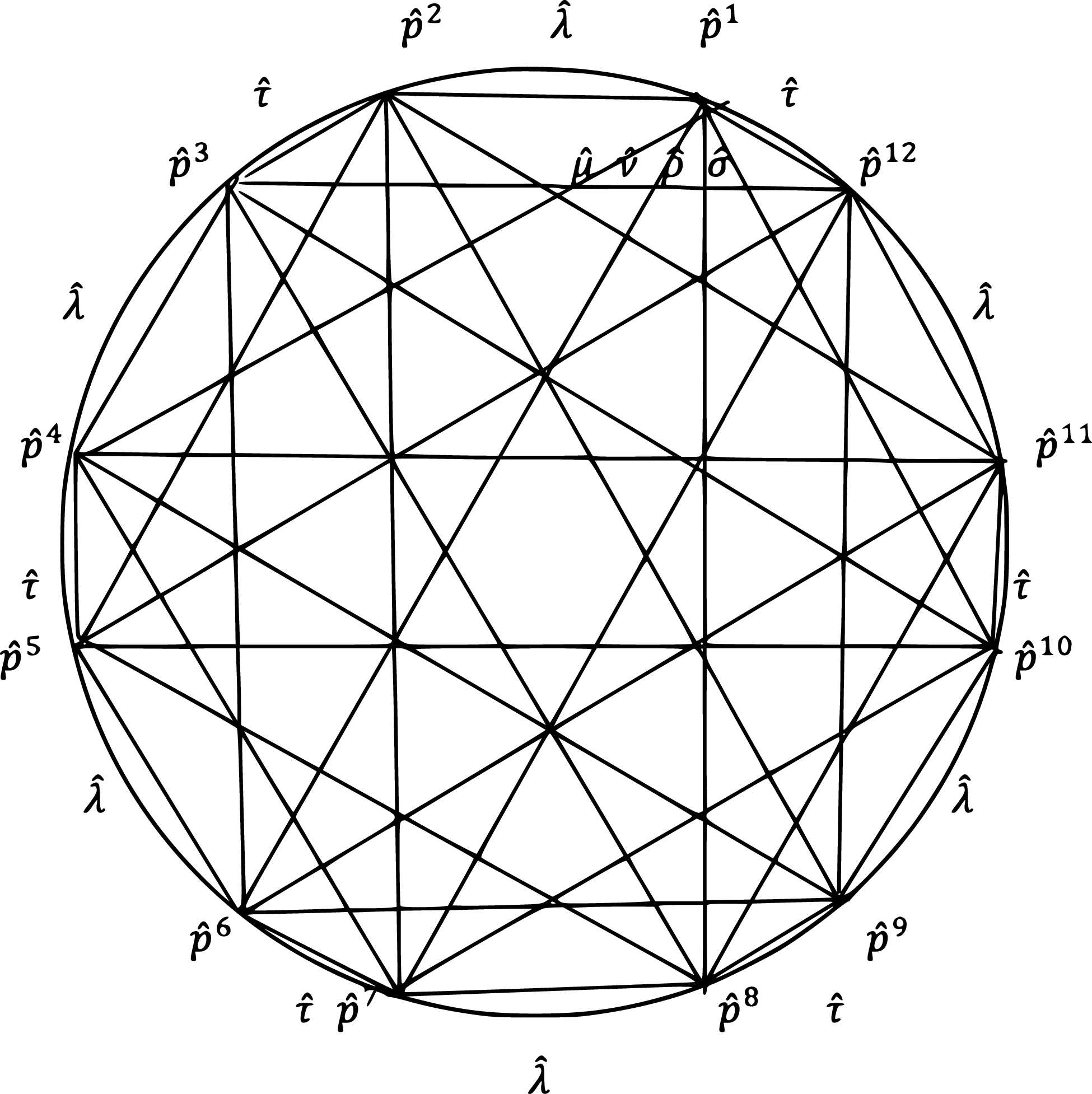}
    \caption{Parallel dodecagon}
 \label{fparag6}
\end{minipage}

\end{figure}

We obtain Proposition \ref{ppi6} under the assumption $\la^1\rho^1=-1=\la^2\rho^2$.  
Since there exists a conformal transformation $\hat C:M\to \hat M$ 
so that the leaves $L_1( \hat p^1)$ and $L_1( \hat p^7)$ are antipodally symmetric, 
we denote this conformal image with hat such as $ \hat M$, $ \hat p^t$, $ \hat\la^t$ etc. 
Then we have $ \hat\la^1 \hat\rho^1=-1= \hat\la^2 \hat\rho^2$ (see Figure \ref{fconf}), and  
can apply the argument in the previous subsection to $ \hat M$, 
since $\hat C$ preserves the Lie curvatures. 
Thus we know $\hat M\cap\hat \ga$ is a parallel dodecagon. 
Let $C:\hat M\to M$ be the inverse conformal transformation of $\hat C$. 

\bep The conformal transformation $C:\hat M\to M$ is an isometry,  
and $M\cap\ga$ is itself a parallel dodecagon. 
\label{ppara6}
\enp
\Proof
Denote $C: \hat M \to M$ restricted to $\R^4_2=\R^2\oplus\R^2_2$ where $\ga\subset \R^2$, by   
$$
C=\bpm *&*&*&0\\
 *&*&*&0\\
 x&y& \check\al&0\\
 0&0&0&1\epm\in O(2,1)\subset O(2,2).
 $$
Apply $C$ to $k_1=\bpm \hat p^t\\1\\0\epm$ and $k_2=\bpm \hat n^t\\0\\1\epm$, and express  
$$
Ck_1=\bpm q^t\\ \hat a_t\\0\epm, \quad Ck_2=\bpm  m^t\\ \hat c_t\\1 \epm. 
$$
Then by (\ref{newpc}), the principal curvatures $\la_i^t$ at the original point $p^t\in M\cap\ga$ are: 
\beq
\la_i^t= \hat a_t\hat\la_i +\hat c_t , \quad i=1,\dots,\quad t=1,\dots 8.  
\label{hatla1}
\enq
The mean curvature $ \hat H^t$ at $\hat p^t\in \hat M\cap\hat\ga$   
dose not depend on $t$ as $\hat M\cap\hat\ga$ is a parallel dodecagon. 
Thus denoting $ \hat H= \hat H^t$, we can express the mean curvature  $H^t$ of $M$ at $p^t$ by  
\beq
H^t= \hat a_t\hat H+ 6m\hat c_t,\quad \hat H=m(\hat \la+\hat\mu+\hat\nu+\hat\rho+\hat\si+\hat\tau), \quad m=1,2, 
\label{CMCH6}
\enq
that is also independent of $t$ as $M$ is CMC. 

We express some vertices of the parallel dodecagon 
$\hat M\cap\hat \ga=\{\hat p^1,\dots,\hat p^{12}\}$ (Figure \ref{fparag6})  
using $\hat u,\hat v>0$ and $\hat k$, $\hat l>0$ s.t.  
\beq
\hat\la^1=\dfrac{\hat v}{\hat u}, \quad  \hat\tau^4=-\dfrac{\hat k}{\hat l},
\label{hatlatau}
\enq 
as 
\beq
\bary\begin{cases}
\hat p^1=\bpm \hat u\\ \hat v\epm=-\hat p^7\\ \hat n^1=\bpm -\hat v\\ \hat u\epm=-\hat n^7\end{cases} &\quad
\begin{cases} \hat p^2=\bpm -\hat u\\ \hat v\epm=\hat p^8\\
\hat n^2=\bpm \hat v\\ \hat u\epm=-\hat n^8,
\end{cases}\\
\begin{cases}
\hat p^4=\bpm -\hat k\\ \hat l\epm=-\hat p^{10}\\
\hat n^4=\bpm \hat l\\ \hat k\epm=-\hat n^{10}\end{cases} &\quad
\begin{cases} \hat p^5=\bpm -\hat k\\-\hat l\epm=-\hat p^{11}\\ \hat n^5=\bpm \hat l\\-\hat k\epm=-\hat n^{11}.  
\end{cases} 
\eary
\label{barpt}
\enq
Thus we have 
\beq
\bary
\begin{cases}
  \hat a_1=  \hat ux+ \hat vy+ \hat\al\\
  \hat c_1=- \hat vx+ \hat u y\end{cases} &\quad 
\begin{cases}
  \hat a_2=- \hat ux+ \hat vy+ \hat\al\\
  \hat c_2= \hat vx+ \hat u y \end{cases} \\
\begin{cases}
 \hat a_4=- \hat kx+\hat ly+ \hat\al\\
 \hat c_4= \hat lx+\hat ky \end{cases} &\quad
 \begin{cases}
 \hat a_5=- \hat kx- \hat ly+ \hat\al\\
 \hat c_5=\hat lx- \hat k y. \end{cases} 
\eary
\label{hatat}
\enq
Since $M$ is CMC,  we have at $p^1$ and $p^2$, 
$$
\bary
0=H^1-H^2 &=(\hat a_1 -\hat a_2)\hat H+6m(\hat c_1-\hat c_2)= 
 2\hat ux( \hat H-6m\hat\la),
\eary
$$
where $\hat u( \hat H-6m\hat\la)<0$  as $\hat\la$ is the largest principal curvature.
Thus  
$$
x=0
$$
follows. Next, at $ p^4$ and $ p^5$  we have
$$
\bary
0&=H^4-H^5
=(\hat a_4-\hat a_5)\hat H+6m(\hat c_4-\hat c_5)\\
&=2\hat ly\hat H+6m\cdot 2\hat ky=2y\hat l( \hat H -6m\hat \tau)
 \eary
 $$
where $\hat l(\hat H-6m\hat\tau)>0$ as $\hat\tau$ is the smallest principal curvature. 
Thus
$$
y=0
$$
holds. 
Therefore,  $C$ is an isometry $\bpm T&0\\0&I'\epm\in   O(2,2)$, $I'=\bpm \eps&0\\0&1\epm$, 
and $M\cap\ga$ itself is a parallel dodecagon $M_\theta\cap\ga_\theta$,  
where $\theta$ is uniquely determined from $H$ of $M$ (Fact 4).  
\qed

\medskip\noindent
{\it Proof of Theorem \ref{tmain} (iv):}
We know that at $p\in M$ where $\la$ takes max or min  and $d_{12}=0$ is satisfied, 
$M\cap\ga$ is isometric to $M_\theta\cap\ga_\theta$ where  
$\theta$ is uniquely determined from the mean curvature $H$ of $M$. 
Thus, $\max\la=\la_\theta=\min\la$ holds, and $\la$ is  constant all over $M$. 
Next, let $p$ be the max or min of $\mu$. As $d_{j1}=0$ and $d_{12}=0$ hold at $p$, 
we can apply the same argument as before to obtain $M\cap\ga$ is again isometric to $M_\theta\cap\ga_\theta$.
 Thus $\max\mu=\mu_\theta=\min\mu$, and $\mu$ is constant all over $M$. 
Then at any point of $M$, $d_{j1}=0$ and $d_{12}=0$ hold, 
and $M\cap\ga$ is isometric to $M_\theta\cap\ga_\theta$ everywhere, 
which implies other principal curvatures $\nu,\rho,\si,\tau$ are constant all over $M$,  
and we proved Theorem \ref{tmain} (iv).

\Remark 
For $\la_1=\cot\theta_1$ where $\theta_1=\dfrac{\pi}{12}+\theta$, $-\dfrac{\pi}{12}<\theta<\dfrac{\pi}{12}$, 
the scalar curvature is given in  \cite{TY} p.147, 
$$
\bary
R_\theta&=36m(m-1)(1+\cot^2(6\theta_1)), 
\eary
$$
which, in the minimal case, becomes
\beq
\bary
R=36m(m-1).
\eary
\label{Rg6}
\enq
When $m=1$ the minimal isoparametric hypersurface is scalar flat, 
and when $m=2$, scalar positive. 

\section{Problems}
\label{sproblem}

We propose: 

\noindent
{\bf Problem 1.} {\em If a closed minimal (or CMC) hypersurface $M$ with $g=4$ 
has constant scalar curvature (CSC), is $M$  Dupin?}

The conclusion is weaker than isoparametric, 
 but if we show this, Chern's conjecture is true for $g=4$ via Theorem \ref{tmain} (ii). 
  
 \noindent
{\bf Problem 2.} {\em If a closed minimal (or CMC) hypersurface $M$ with $g=4$ has constant Lie curvature (CLC), 
is $M$ Dupin?}

If we show this, 
$M$ is isoparametric via Theorem \ref{tmain} (iii). 

\noindent
{\bf Problem 3.} {\em If a closed minimal (or CMC) hypersurface $M$ with $g=6$ has CSC, 
is $M$ Dupin of which Lie curvatures are constant?}

It seems not the case, but if this is true, Chern's conjecture is true for $g=6$ via Theorem \ref{tmain} (iv).  


\medskip
Mathematical Institute, Tohoku University, Sendai 980-8578, Japan

{\it Email address:} r-miyaok\verb|@|tohoku.ac.jp

 \end{document}